\pdfoutput=1
\documentclass{article}

\usepackage{packages}

\newcommand{\cs}[1]{\textcolor{blue}{\small({\sf Clara:} {{#1}})}}

\newcommand{\tail}{\textnormal{t}}
\newcommand{\head}{\textnormal{h}}

\newcommand{\CM}[3]{\textnormal{CM}(\mathcal{H}_{#1}(\vb*{#2}_{#3}))}
\newcommand{\tikzscale}{0.75}
\newcommand{\tikznodesize}{0.5cm}

\title{Hypercurveball algorithm for sampling hypergraphs with fixed degrees}
\author[1]{Yanna J. Kraakman}
\author[1]{Clara Stegehuis}
\affil[1]{University of Twente}
\date{April 7, 2025}

\begin{document}

\maketitle

\begin{abstract}
    Comparative analysis between a network and a random graph model can uncover network properties that significantly deviate from those in random networks. The standard random graph model used for comparison uniformly samples random graphs with the same degrees as the network data, often achieved through edge-swap algorithms. However, for hypergraphs, fewer such methodologies are available. This study introduces the Hypercurveball algorithm, designed to sample random, potentially directed, hypergraphs with fixed degrees. Minor adjustments enable the sampling of hypergraphs without degenerate hyperedges, self-loops, or multi-hyperedges. For most of these algorithms, we prove whether they sample uniformly or with bias. We experimentally show that the Hypercurveball algorithm can be significantly faster or slower than the standard hyperedge-shuffling algorithm, which is the hyperedge-equivalent of the edge-swap algorithm. We present criteria on the hypergraph degree sequence that indicate when the Hypercurveball algorithm is more efficient than the standard hyperedge-shuffling method. Finally, our experimental results suggest polynomial scaling of the mixing time for both the Hypercurveball and hyperedge-shuffling algorithms. 
\end{abstract}

\section{Introduction}
Uniform sampling of graphs with a fixed degree sequence is a fundamental technique to construct null models that can reveal the influence of a given degree sequence on network properties. By contrasting a real-world network to its null model, one can distinguish significant network properties from those attributable to the degree sequence and inherent randomness \cite{newman2018}. A prominent method for such sampling is the edge-rewiring or edge-swapping technique \cite{petersen1891, fosdick2018}, which employs a Markov chain Monte Carlo approach. This method starts with an initial graph and sequentially rewires random pairs of edges for a specified number of iterations. Depending on whether self-loops and multi-edges are allowed, this method can either uniformly generate random graphs or potentially produce certain random graphs with higher probability than others \cite{carstens2017}. The latter is referred to as biased sampling. In scenarios in which uniformity can be reached, the required number of rewirings for achieving uniform randomness, termed the mixing time, has been studied both theoretically for different types of graphs and numerically \cite{levin2009,cooper2007,greenhill2015,greenhill2011,Milo2004,gelman2014}.

More recently, the Curveball method has emerged as a technique for uniform sampling of (di)graphs \cite{verhelst2008,strona2014,Carstens2018}. This method is also a Markov chain Monte Carlo method, and rewires random pairs of \emph{nodes} rather than edges. Experimental evidence suggests that the Curveball method achieves mixing times significantly faster than the edge-swapping method \cite{carstens2015,carstens2016,Carstens2018}.

These sampling methods typically focus on graphs, in which nodes have dyadic interactions. However, many systems exhibit interactions among larger groups of nodes. Therefore, hypergraphs, which generalize the interactions between nodes to larger sets, are often a more accurate representation of a real-world system. They capture the essence of these higher-order interactions better, and often reveal dynamics that are not observable in dyadic graphs \cite{battiston2021,Bick2021,Majhi2022,Lotito2022,Zhang2023}. Hypergraphs can be undirected, such as those modeling social group interactions, or directed, such as those modeling chemical reactions like $A + B \rightarrow C + D$. For both types, a hyperedge-shuffle method, based on the edge-swapping technique, has been proposed and analyzed \cite{chodrow2019,kraakman2024}. Similar to the edge-swapping method for graphs, the hyperedge-shuffle method can either yield uniform random hypergraphs or biased ones, depending on whether self-loops, multi-hyperedges, and/or degenerate hyperedges are allowed during the shuffling process \cite{chodrow2019,kraakman2024}.

In this work, we generalize the Curveball method to accommodate both undirected and directed hypergraphs. Like the Curveball method, this \textit{Hypercurveball} method applies a microcanonical sampling approach, generating hypergraphs whose degree sequence matches the prescribed degree sequence exactly. This approach is stricter than canonical sampling methods, which produce (hyper)graphs where the degree sequence matches the prescribed degree sequence only in expectation. We identify the conditions under which the Hypercurveball method produces uniform random hypergraphs and those where it may introduce bias. Furthermore, we demonstrate that the Hypercurveball algorithm can be more powerful than the existing shuffle algorithm, or vice versa:  for specific degree sequences, the Hypercurveball method yields uniform random samples, whereas the hyperedge-shuffle method results in biased samples, or vice versa. 

For graphs, it has been experimentally demonstrated that the Curveball algorithm mixes significantly faster than the edge swapping method \cite{strona2014,carstens2016}. To our knowledge, there are no theoretical results on the mixing time of the Curveball method for graphs, nor on the mixing time of the hyperedge-shuffle method for hypergraphs. However, there are a number of results on the mixing time of the edge swapping method for graphs. For example, there are polynomial upper bounds for the mixing time of (near-)regular bipartite graphs \cite{kannan1999}, regular graphs \cite{cooper2007}, regular digraphs \cite{greenhill2011} and irregular graphs and digraphs \cite{greenhill2018}. For practical use, however, these bounds seem too large; for instance, the bound for the mixing time of $d$-regular digraphs is of the order $O(d^{26}n^{10}\log(dn))$, where $n$ is the number of nodes. 

We analyze the mixing time experimentally, and show that on real-world network data, often one of the two algorithms is significantly slower. 
We hypothesize that these variations are influenced by a function of the degrees of the nodes and hyperedges in the initial hypergraph, which is supported by our experiments. Moreover, we experimentally show that both methods seem to scale polynomially. Therefore, this work expands the toolkit for hypergraph analysis and offers new insights into the dynamics of hypergraph sampling methods.

\paragraph{Overview of the paper.} Section \ref{section:model} introduces hypergraphs. In Section \ref{section:sampling}, we introduce the Hypercurveball algorithm and show when the method samples uniform random hypergraphs. In Section \ref{section:mixing}, we experimentally compare the mixing time of the Hypercurveball method to the hyperedge-shuffle method and hypothesize that the faster method can be determined by a function of the degrees of the initial hypergraph. The conclusion is in Section \ref{section:conclusion}, and the results are discussed in Section \ref{section:discussion}. Section \ref{section:proofs} contains all proofs.

\section{Hypergraph model}
\label{section:model}
In this section, we introduce undirected and directed hypergraphs, hypergraph spaces and their configuration models. 

\begin{figure}
    \centering
    \newcommand{\offsetx}{1}
\newcommand{\offsety}{0.5}

\begin{tikzpicture}[-Stealth, line width=1.3pt,auto,
                    thick,main node/.style={circle,draw, minimum size=\tikznodesize}, inner sep=1pt]

  \node[main node] (a) at (0,2*\tikzscale) {$a$};
  \node[main node] (b) at (2*\tikzscale,2*\tikzscale) {$b$};
  \node[main node] (c) at (4*\tikzscale,2*\tikzscale) {$c$};
  \node[main node] (d) at (0,0) {$d$};
  \node[main node] (e) at (2*\tikzscale,0) {$e$};
  \node[main node] (f) at (4*\tikzscale,0) {$f$};

\coordinate (adab) at (\tikzscale,\tikzscale);

\draw[-,bend left] (adab) to (a.270);
\draw[-,bend right] (adab) to (d.90);
\draw[out=0, in=0] (adab) to node[right] {$e_1$} (a.0);
\draw[bend right] (adab) to (b.270);

\coordinate (dde) at (\tikzscale,0);

\draw[-,out=180, in=0] (dde) to (d.45);
\draw[-,out=180, in=0] (dde) to (d.315);
\draw (dde) to node[below left] {$e_2$} (e.180);

\draw (b.45) to node[above] {$e_3$} (c.135);
\draw (b.315) to node[below] {$e_4$} (c.225);

\coordinate (cf) at (4.5*\tikzscale,\tikzscale);

\draw[-,out=180, in=270] (cf.180) to (c.270);
\draw[-,out=180, in=90] (cf.180) to (f.90);
\draw[out=0, in=315] (cf) to node[below right] {$e_5$} (c.0);
\draw[out=0, in=45] (cf) to (f.0);

\end{tikzpicture}
    \caption{A directed hypergraph with node set $V=\{a,b,c,d,e,f\}$ and hyperedge set $E=\{e_1,e_2,e_3,e_4,e_5\}$, where $e_1=(\{a,d\},\{a,b\}), e_2=(\{d,d\},\{e\}),e_3=(\{b\},\{c\}), e_4=(\{b\},\{c\})$ and $e_5=(\{c,f\},\{c,f\})$. Hyperedges $e_3$ and $e_4$ are multi-hyperedges, hyperedge $e_2$ is a degenerate hyperedge and hyperedge $e_5$ is a self-loop. The incidence sets of the nodes are $I_a = (\{e_1\},\{e_1\}), I_b = (\{e_3,e_4\},\{e_1\}), I_c = (\{e_5\},\{e_3,e_4,e_5\}), I_d= (\{e_1,e_2,e_2\},\emptyset),I_e = (\emptyset,\{e_2\})$ and $I_f = (\{e_5\},\{e_5\}).$}
    \label{fig:dir_hgraph_example}
\end{figure}
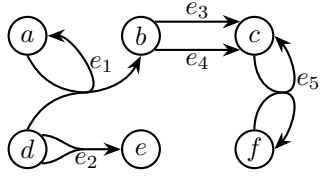

A hypergraph $H=(V,E)$ consist of a node set $V$ and a hyperedge set $E$. Every hyperedge $e\in E$ is a multiset of $V$, that is, it may contain any number of nodes, and each node may appear with any multiplicity. Graphs are a specific type of hypergraph, containing only hyperedges that include exactly two nodes. In undirected hypergraphs, all nodes in a hyperedge have the same role. In directed hypergraphs, each hyperedge $e\in E$ can be written as $(e^{\tail},e^{\head})$, where both $e^{\tail}$ and $e^{\head}$ are multisets of $V$. Thus, each node in a hyperedge is either a tail node (in $e^t$) or a head node (in $e^h$), which can be interpreted as a direction. For example, the chemical reaction $2H_2 + O_2 \rightarrow 2H_2O$ could be modeled as a directed hyperedge with tail nodes $H_2$ (multiplicity 2) and $O_2$ (multiplicity 1), and head node $H_2O$ (multiplicity 2), denoted by $(\{H_2, H_2, O_2\}, \{H_2O,H_2O\})$. Again, digraphs are a specific type of directed hypergraphs, in which each hyperedge contains exactly one tail node and one head node. Note that a hyperedge in a directed hypergraph could contain a node as both a tail and a head node. An example would be a chemical reaction with a catalyst that appears both as a reactant and as a product in the reaction. An illustration of a hypergraph is given in Figure \ref{fig:dir_hgraph_example} \cite{kraakman2024}.

We now define three special hyperedge types that hypergraphs may contain, which we later use to define the configuration models. Hypergraphs can have multi-edges: two hyperedges are considered \emph{multi-hyperedges} if they are equal multisets of nodes (and the nodes have the same roles, in the directed case). Furthermore, a directed hyperedge is classified as a \emph{self-loop} if its multiset of tail nodes equals its multiset of head nodes. For undirected hypergraphs, a self-loop is not defined. While self-loops and multi-edges are also defined for graphs, a hyperedge type that only applies to hypergraphs is the \emph{degenerate hyperedge}: a hyperedge is degenerate if it contains a node with multiplicity of at least two \cite{chodrow2019} (and both have the same role in the hyperedge, in the directed case \cite{kraakman2024}). A degenerate hyperedge can be interpreted as a weighted hyperedge, where each node in the hyperedge has some positive integer weight, as indicated by its multiplicity in the hyperedge. These hyperedge types are illustrated in Figure \ref{fig:dir_hgraph_example}.

Next, we define hypergraph degrees. The degree of a node $v$, denoted by $d_v$, equals the number of times $v$ appears in a hyperedge, i.e., 
\begin{align*}
    d_v := \sum_{e \in E} m_e(v),
\end{align*}
where $m_e(v)$ denotes the multiplicity of node $v$ in hyperedge $e$. The degree of a hyperedge $e$, denoted by $\tilde{d}_e$, equals the size of the multiset that is the hyperedge, i.e., 
\begin{align*}
    \tilde{d}_e := \sum_{v \in V} m_e(v).
\end{align*}
The degree of a hyperedge may also be referred to as the cardinality of the hyperedge. For directed hypergraphs, the multiplicity $m_e(v)$ splits into a tail multiplicity $m_{e^{\tail}}(v)$, counting the multiplicity in $e^t$, and a head multiplicity $m_{e^{\head}}(v)$. The degree of a node or a hyperedge is then a tuple, where the first entry considers only tail multiplicity and the second entry considers only head multiplicity:
\begin{align*}
    d_v := \big(\sum_{e \in E} m_{e^{\tail}}(v), \sum_{e \in E} m_{e^{\head}}(v)\big),\\
    \tilde{d}_e := \big(\sum_{v \in V} m_{e^{\tail}}(v), \sum_{v \in V} m_{e^{\head}}(v)\big).
\end{align*}

Let $\vb*{d}_V := (d_v)_{v \in V}$ and $\vb*{d}_E:= (\tilde{d}_e)_{e \in E}$ denote the sequence of node degrees and hyperedge degrees. The degree sequence $\vb*{d}$ of a hypergraph consists of both of these: $\vb*{d} = (\vb*{d}_V,\vb*{d}_E)$.

To uniformly randomize hypergraph $H$, all hypergraphs with the same degree sequence as $H$ should appear with the same probability. Depending on the problem context, the three special hyperedge types, self-loops ($s$), degenerate hyperedges ($d$) and multi-hyperedges ($m$), may or may not be relevant. For example, a self-loop in a social network might not hold any significance. In that example, one would like to sample uniformly only among the hypergraphs that do not have self-loops. The set of all hypergraphs with some degree sequence $\vb*{d}$ that may contain hyperedge types $x \subseteq \{s,d,m\}$ and not $\{s,d,m\} \backslash x$ is denoted by $\mathcal{H}_x(\vb*{d})$ and is called a hypergraph space.

\begin{definition}[Hypergraph space \cite{kraakman2024}]
    The hypergraph space $\mathcal{H}_{x} (\vb*{d})$, with $x \subseteq \{s,d,m\}$, is a set containing all hypergraphs with degree sequence $\vb*{d}$ that do not contain hyperedge types $s$ (self-loops), $d$ (degenerate hyperedges) or $m$ (multi-hyperedges) that are not in $x$. 
\end{definition}

The most restrictive space for both undirected and directed hypergraphs is $\mathcal{H}(\vb*{d})$, with some degree sequence $\vb*{d}$. The most unrestrictive space is $\mathcal{H}_{d,m}(\vb*{d})$ for undirected hypergraphs and $\mathcal{H}_{s,d,m}(\vb*{d})$ for directed hypergraphs.

A configuration model on hypergraph space $\mathcal{H}_{x} (\vb*{d})$ is defined as the uniform distribution on the hypergraphs in the space $\mathcal{H}_{x} (\vb*{d})$ and is denoted by $\CM{x}{d}{}$. Sampling from $\CM{x}{d}{}$ thus means sampling a hypergraph from the space $\mathcal{H}_{x} (\vb*{d})$ uniformly at random. 

In the remainder of this work, we often make use of incidence sets of nodes to describe a hypergraph. The incidence set of a node $v$, denoted by $I_v$, contains the hyperedges that $v$ participates in (considering multiplicity). More specifically, 
 \begin{align*}
        I_v := \{e: e \in E, v \in e\},
    \end{align*}
    where the multiplicity of the elements is determined by
    \begin{align*}
        m_{I_v}(e) = m_{e}(v).
    \end{align*}

In directed hypergraphs, these incidence sets are again split into a tail incidence set and a head incidence set, i.e., $I_v := (I^{\tail}_v,I^{\head}_v)$. Let $i \in \{\tail,\head\}$, then
    \begin{align*}
        I^{i}_v := \{e: e \in E, v \in e^i\}
    \end{align*}
    and the multiplicity of the elements is determined by 
    \begin{align*}
        m_{I^i_v}(e) = m_{e^i}(v).
    \end{align*}

\begin{figure}[tbp]
    \centering
    \newcommand{\offsetx}{5*\tikzscale}
    
\begin{tikzpicture}[-Stealth, line width=1.3pt,auto,
                    thick,main node/.style={circle,draw, minimum size=\tikznodesize}, inner sep=1pt]

  \node[main node] (a1) at (0,2*\tikzscale) {$a$};
  \node[main node] (b1) at (0,0) {$b$};
  \node[main node] (c1) at (2,\tikzscale) {$c$};

  \node (x) at (\tikzscale, -1*\tikzscale) {$H_1$};
  \node[align=left] at (1.5*\tikzscale, -3*\tikzscale) {$e_1 = (\{b\},\{a\})$\\$e_2=(\{a,b\},\{c\})$\\ $I_a = (\{e_2\},\{e_1\})$ \\ $I_b = (\{e_1,e_2\},\emptyset)$\\$I_c = (\emptyset, \{e_2\})$};

\draw (b1) to node {$e_1$} (a1) ;

\coordinate (abc1) at (\tikzscale,\tikzscale); 

\draw[-, bend right] (a1.300) to (abc1.180);
\draw[-, bend left] (b1.60) to (abc1.180);
\draw[] (abc1) to (c1);

\node[above] at (abc1) {$e_2$};

  \node[main node] (a2) at (\offsetx,2*\tikzscale) {$a$};
  \node[main node] (b2) at (\offsetx,0) {$b$};
  \node[main node] (c2) at (\offsetx+2,\tikzscale) {$c$};

  \node (x) at (\offsetx + \tikzscale, -1*\tikzscale) {$H_2$};
  \node[align=left] at (\offsetx + 1.5*\tikzscale, -3*\tikzscale) {$e'_1 = (\{a\},\{c\})$\\$e'_2=(\{b,b\},\{a\})$\\ $I_a = (\{e'_1\},\{e'_2\})$ \\ $I_b = (\{e'_2,e'_2\},\emptyset)$\\$I_c = (\emptyset, \{e'_1\})$};

\draw (a2) to node {$e'_1$} (c2) ;

\coordinate (bba2) at (\offsetx,\tikzscale); 

\draw[-, bend right, in = 180] (b2.45) to (bba2.270);
\draw[-, bend left, in=180] (b2.135) to (bba2.270);
\draw[] (bba2) to (a2);

\node[left, xshift=-5*\tikzscale pt] at (bba2) {$e'_2$};

\end{tikzpicture}
    \caption{Changing the incidence sets of a hypergraph also changes the edge set.}
    \label{fig:E_and_Inc}
\end{figure}
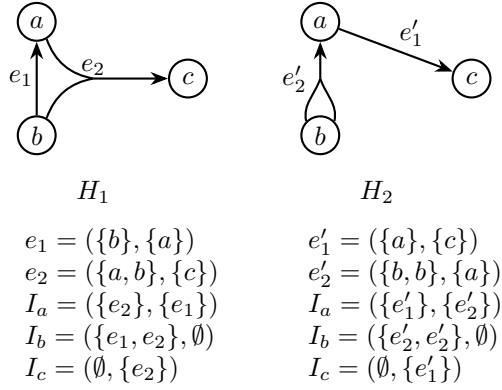

The hyperedges are multisets of nodes. Therefore, an incidence set is a multiset of multisets of nodes. The edge set $E$ and the set of incidence sets $\{I_v: v \in V\}$ are different ways to describe the same hypergraph. Therefore, defining one also fixes the other. In particular, changing the edge set also changes the incidence sets, and vice versa. For example, in $H_1$ in Figure \ref{fig:E_and_Inc} the incidence set of node $a$ is $I_a(H_1) = (\{e_2\},\{e_1\})=(\{(\{a,b\},\{c\})\},\{(\{b\},\{a\})\})$. When changing the hypergraph such that node $a$ becomes incident to $e_1$ (tail) and $e_2$ (head), and adjusting the incidences of nodes $b$ and $c$ to keep the degrees of $e_1$ and $e_2$ fixed (hypergraph $H_2$ in Figure \ref{fig:E_and_Inc}), we would like to write $I_a(H_2) = (\{e_1\},\{e_2\})$. However, by changing the incidences, $e_1$ and $e_2$ change into $e'_1$ and $e'_2$, such that $I_a(H_2) = (\{e'_1\},\{e'_2\}) = (\{(\{a\},\{c\})\},\{(\{b,b\},\{a\})\})$.

The interplay between the edge set and the incidence sets representations can make the interpretation of the incidence sets difficult. Therefore, we choose the incidence set notation as leading and assume that the edge set notation will automatically be adjusted to follow the changes in the incidence set. For example, in Figure~\ref{fig:E_and_Inc}, $I_a(H_2) = (\{e_1\},\{e_2\})$ describes the transformation of $H_1$ into $H_2$: node $a$ becomes adjacent to the edge with label $e_1$ (tail) and the edge with label $e_2$ (head). The new hyperedges behind the labels $e_1$ and $e_2$ are implicit in this notation, and can be derived by looking at all incidence sets if needed. However, often the exact form of the hyperedges is not relevant. In the remainder of this work, for the sake of readability, the elements of the incidence sets should be interpreted as the labels of the incident hyperedges. These labels do not have any meaning other than to distinguish different hyperedges, so that exchanging the labels of two hyperedges does not change the hypergraph. The explicit dependency of the incidence sets on the hypergraph is omitted if it is clear from the context.
    
For the experimental part of this work, we implemented hypergraphs as incidence sets. In particular, a hypergraph is described by a dictionary, where each node is a key whose value is the incidence set of this node. The implementation can be found here: \cite{algorithms}. In addition, we provide the implementation of a hypergraph as hyperedge sets, as we use this to experimentally analyze the hyperedge-shuffle method. One representation can easily be converted into the other. A different way to implement a hypergraph could be using an incidence matrix, with a row for each node and a column for each hyperedge. Each entry in the matrix indicates the multiplicity of a node in a hyperedge, and if the hypergraph is directed then the entries in the matrix are tuples. This implementation may be more intuitive than the incidence sets, but is computationally more expensive due to the large size of the matrix. The incidence matrix could be implemented as a sparse matrix, but this yields computational difficulties when frequently updating values in the matrix, as we do in the remainder of this work.

\section{Sampling: Hypercurveball algorithm}
\label{section:sampling}
Section \ref{section:curveball} introduces our Hypercurveball algorithms for hypergraphs. In Section \ref{section:uniformity}, we show when these algorithms sample uniformly from the hypergraph spaces and when they may sample with bias. We compare these results to similar results for the hyperedge-shuffle algorithm in Section \ref{section:comparison}.

\subsection{Hypercurveball algorithm}
\label{section:curveball}
We introduce the Hypercurveball algorithm for hypergraphs, which generalizes the Curveball algorithm for graphs and digraphs \cite{carstens2016}. The essence of the Curveball algorithm for graphs is to randomly pick two nodes and reconfigure the edges incident to these nodes while keeping their degrees fixed. This process, known as a \emph{trade}, contrasts with the edge-swapping algorithm where two edges are randomly selected, and their incident nodes are reconfigured \cite{bollobas1980}.

The Curveball algorithm has experimentally been shown to mix much faster than the edge-swapping algorithm \cite{carstens2015,carstens2016,Carstens2018}. This is explained by the fact that an edge swap only changes four incidences (two edges with two incidences each), whereas a Curveball trade can change multiple incidences at once, since a node can be incident to many edges. Therefore, the Curveball algorithm often randomizes an initial graph more quickly.

We now introduce a hypertrade for hypergraphs. Two nodes $v$ and $w$ are picked uniformly at random, and their incidence sets $I_v(H)$ and $I_w(H)$ are randomly reconfigured such that the degrees of the hypergraph remain fixed. Let $A \uplus B$ denote the multiset union with multiplicity $m_{A \uplus B}(x)=m_A(x)+m_B(x)$. For example, if $A=\{x,y\}$ and $B=\{y,z\}$ then $A \uplus B=\{x,y,y,z\}$. $I_v(H)$ and $I_w(H)$ are reconfigured by creating a new multiset $I_{v+w} := I_v(H) \uplus I_w(H)$ and randomly repartitioning it into two sets $I'_v:=I_v(H')$, $I'_w:=I_w(H')$ of sizes $|I_v|$ and $|I_w|$, respectively, where $H'$ is the hypergraph described by the new incidence sets. By creating the partition, the degrees of all hyperedges stay fixed. By keeping the sizes of the incidence sets fixed, the number of hyperedges that $v$ resp. $w$ participates in remains fixed. Therefore, after applying any number of hypertrades, the resulting hypergraph has the same degrees as the initial hypergraph.

\begin{figure}[tbp]
    \centering
    \newcommand{\offsetx}{5*\tikzscale}
    
\begin{tikzpicture}[-Stealth, line width=1.3pt,auto,
                    thick,main node/.style={circle,draw, minimum size=\tikznodesize}, inner sep=1pt]

  \node[main node] (a1) at (0,4*\tikzscale) {$a$};
  \node[main node] (b1) at (3*\tikzscale,4*\tikzscale) {$b$};
  \node[main node] (c1) at (0,2*\tikzscale) {$c$};
  \node[main node] (d1) at (1.5*\tikzscale,2*\tikzscale) {$d$};
  \node[main node] (e1) at (3*\tikzscale,2*\tikzscale) {$e$};
  \node[main node] (f1) at (1.5*\tikzscale,0) {$f$};

  \node (x) at (1.5*\tikzscale, -\tikzscale) {$H_1$};
  \node[align=left] at (1.5*\tikzscale, -2*\tikzscale) {$I_c = (\{e_3,e_3\},\{e_2\})$ \\ $I_e = (\{e_2\},\{e_1\})$};

\coordinate (abe1) at (1.5*\tikzscale,3*\tikzscale); 

\draw[-, bend left] (abe1) to (a1.315);
\draw[bend right] (abe1) to (b1.225);
\draw[bend left] (abe1) to (e1.90);

\node[above] at (abe1) {$e_1$};

\coordinate (ecd1) at (2.25*\tikzscale,2.5*\tikzscale); 

\draw[-] (ecd1) to (e1.135);
\draw[bend right] (ecd1) to (c1.90);
\draw[bend right] (ecd1) to (d1.90);

\node[below] at (ecd1) {$e_2$};

\coordinate (ccf1) at (0.75*\tikzscale,\tikzscale); 

\draw[-,bend right, out=0] (ccf1) to (c1.325);
\draw[-,bend left, out=0] (ccf1) to (c1.270);
\draw (ccf1) to (f1.135);

\node[right] at (ccf1) {$e_3$};

 \node[main node] (a2) at (\offsetx,4*\tikzscale) {$a$};
  \node[main node] (b2) at (\offsetx+3*\tikzscale,4*\tikzscale) {$b$};
  \node[main node] (c2) at (\offsetx,2*\tikzscale) {$c$};
  \node[main node] (d2) at (\offsetx+1.5*\tikzscale,2*\tikzscale) {$d$};
  \node[main node] (e2) at (\offsetx+3*\tikzscale,2*\tikzscale) {$e$};
  \node[main node] (f2) at (\offsetx+1.5*\tikzscale,0) {$f$};

  \node (x) at (\offsetx + 1.5*\tikzscale, -\tikzscale) {$H_2$};
\node[align=left] at (\offsetx + 1.5*\tikzscale, -2*\tikzscale) {$I_c = (\{e_2,e_3\},\{e_2\})$ \\ $I_e = (\{e_3\},\{e_1\})$};

\coordinate (abe2) at (\offsetx+1.5*\tikzscale,3*\tikzscale); 

\draw[-, bend left] (abe2) to (a2.315);
\draw[bend right] (abe2) to (b2.225);
\draw[bend left] (abe2) to (e2.90);

\node[above] at (abe2) {$e_1$};

\coordinate (ccd2) at (\offsetx+ 0.75*\tikzscale,2.5*\tikzscale); 

\draw[-, bend left] (ccd2) to (c2.0);
\draw[bend right, out=270] (ccd2) to (c2.90);
\draw[bend left, out=90] (ccd2) to (d2.90);

\node[above, yshift = 4.25 * \tikzscale pt] at (ccd2) {$e_2$};

\coordinate (cef2) at (\offsetx + 1.5*\tikzscale,0.75*\tikzscale); 

\draw[-,bend left, out=300] (cef2) to (c2.270);
\draw[-,bend right, out=60] (cef2) to (e2.270);
\draw (cef2) to (f2.90);

\node[above, yshift = 9*\tikzscale pt] at (cef2) {$e_3$};

\node[main node] (a3) at (2*\offsetx,4*\tikzscale) {$a$};
  \node[main node] (b3) at (2*\offsetx+3*\tikzscale,4*\tikzscale) {$b$};
  \node[main node] (c3) at (2*\offsetx,2*\tikzscale) {$c$};
  \node[main node] (d3) at (2*\offsetx+1.5*\tikzscale,2*\tikzscale) {$d$};
  \node[main node] (e3) at (2*\offsetx+3*\tikzscale,2*\tikzscale) {$e$};
  \node[main node] (f3) at (2*\offsetx+1.5*\tikzscale,0) {$f$};

  \node (x) at (2*\offsetx + 1.5*\tikzscale, -\tikzscale) {$H_3$};

\node[align=left] at (2*\offsetx + 1.5*\tikzscale, -2*\tikzscale) {$I_c = (\{e_3,e_3\},\{e_1\})$ \\ $I_e = (\{e_2\},\{e_2\})$};

\coordinate (abc3) at (2*\offsetx+0.75*\tikzscale,3*\tikzscale); 

\draw[-, bend left] (abc3) to (a3.315);
\draw[bend right] (abc3) to (b3.225);
\draw[bend right, out=90, in=140] (abc3) to (c3.45);

\node[above right] at (abc3) {$e_1$};

\coordinate (ded3) at (2*\offsetx+ 2.25*\tikzscale,2.5*\tikzscale); 

\draw[-, bend right] (ded3) to (e3.180);
\draw[bend left, out=90] (ded3) to (e3.90);
\draw[bend right, out=270] (ded3) to (d3.90);

\node[above, yshift = 4.25*\tikzscale pt] at (ded3) {$e_2$};

\coordinate (ccf3) at (2*\offsetx + 0.75*\tikzscale,\tikzscale); 

\draw[-,bend right, out=0] (ccf3) to (c3.325);
\draw[-,bend left, out=0] (ccf3) to (c3.270);
\draw (ccf3) to (f3.135);

\node[right] at (ccf3) {$e_3$};

\node[main node] (a4) at (3*\offsetx,4*\tikzscale) {$a$};
  \node[main node] (b4) at (3*\offsetx+3*\tikzscale,4*\tikzscale) {$b$};
  \node[main node] (c4) at (3*\offsetx,2*\tikzscale) {$c$};
  \node[main node] (d4) at (3*\offsetx+1.5*\tikzscale,2*\tikzscale) {$d$};
  \node[main node] (e4) at (3*\offsetx+3*\tikzscale,2*\tikzscale) {$e$};
  \node[main node] (f4) at (3*\offsetx+1.5*\tikzscale,0) {$f$};

  \node (x) at (3*\offsetx + 1.5*\tikzscale, -\tikzscale) {$H_4$};
  \node[align=left] at (3*\offsetx + 1.5*\tikzscale, -2*\tikzscale) {$I_c = (\{e_2,e_3\},\{e_1\})$ \\ $I_e = (\{e_3\},\{e_2\})$};

\coordinate (abc4) at (3*\offsetx+0.75*\tikzscale,3.5*\tikzscale); 

\draw[-] (abc4) to (a4.315);
\draw[bend right] (abc4) to (b4.225);
\draw[bend right, out=90, in=180] (abc4) to (c4.90);

\node[above right] at (abc4) {$e_1$};

\coordinate (cde4) at (3*\offsetx+ 0.75*\tikzscale,2.5*\tikzscale); 

\draw[-] (cde4) to (c4.45);
\draw[bend left] (cde4) to (d4.90);
\draw[bend left] (cde4) to (e4.90);

\node[below, yshift = -1.5*\tikzscale pt] at (cde4) {$e_2$};

\coordinate (cef4) at (3*\offsetx + 1.5*\tikzscale,0.75*\tikzscale); 

\draw[-,bend left, out=300] (cef4) to (c4.270);
\draw[-,bend right, out=60] (cef4) to (e4.270);
\draw (cef4) to (f4.90);

\node[above, yshift = 9*\tikzscale pt] at (cef4) {$e_3$};

\end{tikzpicture}
    \caption{Performing \texttt{hypertrade(}$c,e$\texttt{)} on $H_1$ can result in the hypergraph $H_1,H_2,H_3$ or $H_4$.}
    \label{fig:example_trade}
\end{figure}
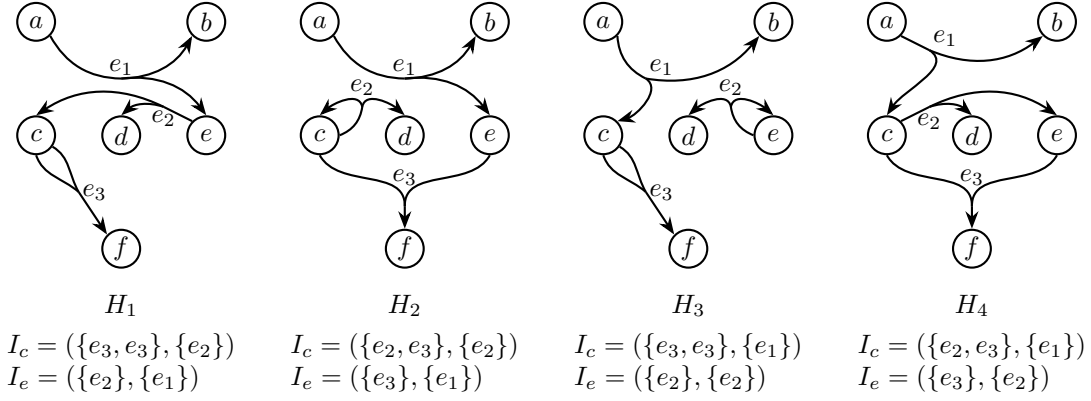

For directed hypergraphs, a hypertrade consists of first reconfiguring $I^{\tail}_{v}(H)$ and $I^{\tail}_{w}(H)$, followed by reconfiguring $I^{\head}_{v}(H)$ and $I^{\head}_{w}(H)$. This ensures that nodes do not switch roles from tail to head, or vice versa. Algorithm \ref{alg:trade_d} provides an implementation of the hypertrade. A hypertrade on a directed hypergraph is illustrated in Figure \ref{fig:example_trade}.

\begin{algorithm}[tb]
\SetKwFunction{Fun}{hypertrade}
\SetKwProg{Fn}{Function}{:}{}
\SetKwInOut{Input}{Input}
\SetKwInOut{Output}{Output}
\Fn{\Fun{v,w}}{
    \Input{$v,w,I_v,I_w$}
    \Output{$I'_v,I'_w$}
    \eIf{\textnormal{undirected hypergraph}}{
    toReconfigure := $\{(I_v,I_w)\}$
    }{
    toReconfigure := $\{(I_v^{\tail},I_w^{\tail}),(I_v^{\head},I_w^{\head})\}$
    }
    \For{$(I,J)$ in \textnormal{toReconfigure}}{
    Create random partition of $I \uplus J$ into two sets $I'$ and $J'$ of sizes $|I|$ and $|J|$, resp. \label{line:partition}
    }
    \Return{$I'_v,I'_w$}
}
 \caption{hypertrade allowing degenerate hyperedges, self-loops and multi-hyperedges}
 \label{alg:trade_d}
\end{algorithm}

Algorithm \ref{alg:curveball_sdm} shows the full Hypercurveball algorithm. The algorithm performs $N$ random trades. To fully randomize an initial hypergraph $H$, $N$ should be at least the mixing time of the Hypercurveball algorithm, which most likely depends on the number of nodes, the number of edges and the degrees of $H$, but is not known theoretically. In practice, for similar Markov chain methods on graphs, often $N=100 \cdot |E|$ is used \cite{Milo2004}. For hypergraphs, to the best of our knowledge, there are no theoretical nor experimental results on the mixing time of such Markov chain Monte Carlo methods. In Section \ref{section:mixing}, we present some experimental results on the mixing time for hypergraphs. 

\begin{algorithm}[tb]
\textbf{Input:} hypergraph $H=(V,E)$, large number $N \in \mathbb{N}$ \label{line:N}\\
 \For{$i$ in $1,\hdots,N$}{
 Select two nodes $v,w \in V$ uniformly at random\\
 Let $I'_v,I'_w$ = \texttt{hypertrade}($v,w)$ \label{line:trade}\\
 $I_v \gets I'_v$\\
 $I_w \gets I'_w$
 }
 \Return{H}
 \caption{Hypercurveball algorithm for $\mathcal{H}_{s,d,m} (\vb*{d})$}
 \label{alg:curveball_sdm}
\end{algorithm}

Algorithm \ref{alg:curveball_sdm} samples from $\CM{s,d,m}{d}{}$, the hypergraph space that allows self-loops, degenerate hyperedges and multi-hyperedges. To sample from a hypergraph space without degenerate hyperedges, \texttt{hypertrade} can be modified to avoid creating degenerate hyperedges. This is achieved by altering line \ref{line:partition} in Algorithm \ref{alg:trade_d} such that the sets $I'$ and $J'$ do not contain duplicate elements. It keeps the elements that are both in $I$ and $J$ fixed and only reconfigures the elements in the symmetric difference $(I \backslash J) \cup (J \backslash I)$. The result is \texttt{hypertrade\_nodeg}, as presented in Algorithm \ref{alg:trade_no_d}.

\begin{algorithm}[tbp]
\SetKwFunction{Fun}{hypertrade\_nodeg}
\SetKwProg{Fn}{Function}{:}{}
\SetKwInOut{Input}{Input}
\SetKwInOut{Output}{Output}
\Fn{\Fun{v,w}}{
    \Input{$v,w,I_v,I_w$}
    \Output{$I'_v,I'_w$}
    \eIf{\textnormal{undirected hypergraph}}{
    toReconfigure := $\{(I_v,I_w)\}$
    }{
    toReconfigure := $\{(I_v^{\tail},I_w^{\tail}),(I_v^{\head},I_w^{\head})\}$
    }
    \For{$(I,J)$ in \textnormal{toReconfigure}}{
    Let $\Tilde{I} = I \backslash J$\\
    Let $\Tilde{J} = J \backslash I$\\
    Create random partition of $\Tilde{I} \cup \Tilde{J}$ into two sets $\Tilde{I}'$ and $\Tilde{J}'$ of sizes $|\Tilde{I}|$ and $|\Tilde{J}|$, resp.\\
    Let $I' = \Tilde{I}' \cup (I \cap J)$\\
    Let $J' = \Tilde{J}' \cup (I \cap J)$
    }
    \Return{$I'_v,I'_w$}
}
 \caption{hypertrade not allowing degenerate hyperedges}
 \label{alg:trade_no_d}
\end{algorithm}

\texttt{hypertrade\_nodeg} can also be adjusted to sample from $\mathcal{H}_{s}(\vb*{d})$, an undirected hypergraph space without degenerate hyperedges or multi-hyperedges. To do this, we first identify hyperedge pairs $\{x,y\}$ which could become multi-hyperedges after a trade involving nodes $v$ and $w$. These are the hyperedges $x,y$ for which
\begin{align}
\label{potential_multihyperedges}
    &x,y \in I_v \cup I_w\\
    &\forall u \in V \backslash \{v,w\}: x \in I_u \iff y \in I_u. \nonumber
\end{align}
Then $x$ and $y$ should not be reconfigured into the same incidence set $I'_v$ or $I'_w$. To ensure this, for each such pair $\{x,y\}$, one hyperedge is assigned to $I'_v$ and one hyperedge is assigned to $I'_w$, uniformly at random. A hyperedge $x$ cannot appear in multiple such pairs, since that would imply that its paired hyperedges are multi-hyperedges. Afterwards, the remaining hyperedges in $I_v$ and $I_w$ are reconfigured. This results in \texttt{hypertrade\_simple}, see Algorithm \ref{alg:trade_no_d_no_m}.  

\begin{algorithm}[tb]
\SetKwFunction{Fun}{hypertrade\_simple}
\SetKwProg{Fn}{Function}{:}{}
\SetKwInOut{Input}{Input}
\SetKwInOut{Output}{Output}
\Fn{\Fun{v,w}}{
    \Input{$v,w,\forall u\in V: I_u$}
    \Output{$I'_v,I'_w$}
    Let $T$ be the set of pairs of hyperedges satisfying (\ref{potential_multihyperedges})\\
    $I'_v \gets \emptyset$\\
    $I'_w \gets \emptyset$\\
    \For{$\{x,y\}$ in T}{
    $\begin{cases}
        I'_v \gets I'_v \cup \{x\}, I'_w \gets I'_w \cup \{y\}& \textnormal{w.p. }1/2\\
        I'_v \gets I'_v \cup \{y\}, I'_w \gets I'_w \cup \{x\} & \textnormal{w.p. }1/2
    \end{cases}$}
    Let $\Tilde{I}_v = I_v \backslash (I_w \cup \{x: \{x,y\} \in T\})$\\
    Let $\Tilde{I}_w = I_w \backslash (I_v \cup \{x: \{x,y\} \in T\})$\\
    Create random partition of $\Tilde{I}_v \cup \Tilde{I}_w$ into two sets $\Tilde{I}'_v$ and $\Tilde{I}'_w$ of sizes $|\Tilde{I}_v|$ and $|\Tilde{I}_v|$, resp.\\
    $I'_v \gets I'_v \cup \Tilde{I}'_v \cup (I_v \cap I_w)$\\
    $I'_w \gets I'_w \cup \Tilde{I}'_w \cup (I_v \cap I_w)$\\
    \Return{$I'_v,I'_w$}
}
 \caption{hypertrade for undirected hypergraphs not allowing degenerate hyperedges nor multi-hyperedges}
 \label{alg:trade_no_d_no_m}
\end{algorithm}

The presented \texttt{hypertrade}, \texttt{hypertrade\_nodeg} and \texttt{hypertrade\_simple} are unbiased. That is, they ouput all possible incidence sets $I'_v$ and $I'_w$ with equal probability. This can easily be checked. To sample from the remaining hypergraph spaces ($\mathcal{H}_{s,d}, \mathcal{H}_{d,m}, \mathcal{H}_s$ (directed), $\mathcal{H}_d, \mathcal{H}_m$ and $\mathcal{H}$), the \texttt{hypertrade} function cannot easily be adjusted to output incidence sets without creating bias (see Appendix \ref{app:bias}). Instead, an additional line is inserted after line \ref{line:trade} in Algorithm \ref{alg:curveball_sdm} to check for self-loops and/or multi-hyperedges. If none are found, the hypertrade is accepted. If the hypertrade is not accepted, it still counts towards the total of $N$ trades. The adjusted algorithm for the space without self-loops is shown in Algorithm \ref{alg:curveball_overig}. Similar adjustments can be made for the other spaces.

\begin{algorithm}[tb]
\textbf{Input:} hypergraph $H=(V,E)$, large number $N \in \mathbb{N}$\\
 \For{$i$ in $1,\hdots,N$}{
 Select two nodes $v,w \in V$ uniformly at random\\
 \eIf{\textnormal{degenerate allowed}}{
 Let $I'_v,I'_w$ = \texttt{hypertrade}($v,w)$
 }{
 Let $I'_v,I'_w$ = \texttt{hypertrade\_nodeg}($v,w)$
 }
 \If{\textnormal{Resulting hypergraph with} $I_v = I'_v$ \textnormal{and} $I_w = I'_w$ \label{line:check}\textnormal{contains no self-loops}}{
 $I_v \gets I'_v$\\
 $I_w \gets I'_w$
 } 
 }
 \Return{H}
 \caption{Hypercurveball algorithm for $\mathcal{H}_{d,m} (\vb*{d})$ and $\mathcal{H}_{m} (\vb*{d})$}
 \label{alg:curveball_overig}
\end{algorithm}

In summary, the following Hypercurveball algorithms apply to each hypergraph space, for any undirected degree sequence $\vb*{d}_1$ or directed degree sequence $\vb*{d}_2$:
\begin{itemize}
    \item Undirected space $\mathcal{H}(\vb*{d}_1)$: Algorithm \ref{alg:curveball_sdm} using \texttt{hypertrade\_simple};
    \item Undirected space $\mathcal{H}_{d}(\vb*{d}_1)$ and directed spaces $\mathcal{H}_s(\vb*{d}_2)$ and $\mathcal{H}_{s,d}(\vb*{d}_2)$: Algorithm \ref{alg:curveball_overig} with line 8 replaced by a check for multi-hyperedges;
    \item Undirected space $\mathcal{H}_m(\vb*{d}_1)$ and directed space $\mathcal{H}_{s,m}(\vb*{d}_2)$: Algorithm \ref{alg:curveball_sdm} using \texttt{hypertrade\_nodeg};
    \item Undirected space $\mathcal{H}_{d,m}(\vb*{d}_1)$ and directed space $\mathcal{H}_{s,d,m}(\vb*{d}_2)$: Algorithm \ref{alg:curveball_sdm} using \texttt{hypertrade};
    \item Directed spaces $\mathcal{H}(\vb*{d}_2)$ and $\mathcal{H}_d(\vb*{d}_2)$: Algorithm \ref{alg:curveball_overig} with extra check for multi-hyperedges in line \ref{line:check};
    \item Directed spaces $\mathcal{H}_{m}(\vb*{d}_2)$ and $\mathcal{H}_{d,m}(\vb*{d}_2)$: Algorithm \ref{alg:curveball_overig}.
\end{itemize}

The implemented Hypercurveball algorithm can be found here:  \cite{algorithms}.

\subsection{Uniform sampling}
\label{section:uniformity}

The algorithms discussed in Section \ref{section:curveball} effectively sample from the specified hypergraph spaces. However, it is not straightforward whether these samples are uniformly distributed, meaning whether they are drawn from the associated configuration models. In this section, we identify the hypergraph spaces where uniform sampling is guaranteed. All supporting proofs are provided in Section \ref{section:proofs}. Initially, we determine the hypergraph spaces for which the Hypercurveball method produces uniform samples across all degree sequences of the initial hypergraph.
\\
\begin{theorem}
    \label{thm:uniform_sampling}
    The Hypercurveball algorithm samples uniformly from
    \begin{itemize}
        \item the undirected hypergraph spaces $\mathcal{H}_{d,m}(\vb*{d}_1), \mathcal{H}_{m}(\vb*{d}_1)$,
        \item the directed hypergraph spaces $\mathcal{H}_{s,d,m}(\vb*{d}_2), \mathcal{H}_{s,m}(\vb*{d}_2)$,
    \end{itemize}
    for any undirected hypergraph degree sequence $\vb*{d}_1$ and directed hypergraph degree sequence $\vb*{d}_2$.
\end{theorem}

Secondly, we identify the hypergraph spaces where the Hypercurveball method consistently achieves uniform sampling when the initial hypergraph has hyperedges of size 2, i.e., when it is a (di)graph. Note that Theorem \ref{thm:uniform_sampling} holds for (di)graph as well, so the Curveball algorithm samples uniformly from the undirected hypergraph spaces $\mathcal{H}_{d,m}(\vb*{d}'_1)$ (where a degenerate edge is known as a self-loop) and $\mathcal{H}_{m}(\vb*{d}'_1)$ and the directed hypergraph space $\mathcal{H}_{s,m}(\vb*{d}'_2)$. The additional results for graph spaces contain the proof for the undirected space $\mathcal{H}(\vb*{d}'_1)$ by \cite{carstens2016} and the proof for the directed space $\mathcal{H}_{s}(\vb*{d}'_2)$ by \cite{carstens2015}. We extend these proofs to the settings without multi-edges and without multi-edges, self-loops and degenerate edges. This way, all cases for (di)graphs are solved.
\\

\begin{theorem}
    \label{thm:uniform_sampling_graph}
    The Curveball algorithm samples uniformly from
    \begin{itemize}
        \item all undirected graph spaces $\mathcal{H}_x(\vb*{d}'_1)$ except the space with self-loops and without multi-edges ($\mathcal{H}_d(\vb*{d}'_1)$),
        \item all directed graph spaces $\mathcal{H}_x(\vb*{d}'_2)$,
    \end{itemize}
    for any $x \subseteq \{s,d,m\}$, any graph degree sequence $\vb*{d}'_1$ and digraph degree sequence $\vb*{d}'_2$.
\end{theorem}

Third, we identify the hypergraph spaces where the Hypercurveball method may introduce sampling bias, depending on the degree sequence of the initial hypergraph.
\\
\begin{theorem}
    \label{thm:biased_sampling}
    There exist degree sequences $\tilde{\vb*{d}}, \hat{\vb*{d}}$ and $\vb*{d}^*$ for which the Hypercurveball algorithm samples with bias from 
    \begin{itemize}
        \item the undirected hypergraph space $\mathcal{H}_d(\tilde{\vb*{d}})$,
        \item the directed hypergraph spaces $\mathcal{H}_{s,d}(\hat{\vb*{d}}), \mathcal{H}_{d,m}(\vb*{d}^*),\mathcal{H}_d(\vb*{d}^*)$.
    \end{itemize}
\end{theorem}

The results are summarized in Table \ref{tab:results_undir} for undirected hypergraphs and in Table \ref{tab:results_dir} for directed hypergraphs. Here `No' in the `Always uniform?' column implies that at least one degree sequence exists for which the Hypercurveball algorithm samples with bias. As can be seen in the tables, there are a few open cases: the undirected space $\mathcal{H}(\vb*{d}_1)$ and the directed spaces $\mathcal{H}(\vb*{d}_2), \mathcal{H}_m(\vb*{d}_2)$ and $\mathcal{H}_s(\vb*{d}_2)$. For these cases, we generated all possible degree sequences up to some maximal size and found that the Hypercurveball algorithm sampled from each such degree sequence uniformly. The maximal size that we analyzed is shown in the notes in the tables. Most examples of degree sequences that yield biased sampling are rather small, hence it is likely that the open cases can all be sampled uniformly using the Hypercurveball method.

Notably, the results for the directed hypergraph spaces $\mathcal{H}_d(\vb*{d}_2), \mathcal{H}_{d,m}(\vb*{d}_2),\mathcal{H}_{s,d}(\vb*{d}_2),\mathcal{H}_{s,m}(\vb*{d}_2)$ and $\mathcal{H}_{s,d,m}(\vb*{d}_2)$ align with those of the hyperedge-shuffle method for directed hypergraphs \cite{kraakman2024}. However, for at least $\mathcal{H}_d(\vb*{d}_2)$ and $\mathcal{H}_{d,m}(\vb*{d}_2)$, the two methods are not equivalent in terms of uniform sampling: given a degree sequence and a directed hypergraph space, one method may sample uniformly while the other does not.\\
\begin{theorem}
    \label{thm:methods_not_equivalent}
    There exists a degree sequence $\vb*{d}^*$ for which the Hypercurveball method samples uniformly and the hyperedge-shuffle method samples with bias from $\mathcal{H}_{d}(\vb*{d}^*)$ and $\mathcal{H}_{d,m} (\vb*{d}^*_2)$, and there exists a degree sequence $\vb*{d}^{**}$ for which the reverse is true.
\end{theorem}

The case of $\mathcal{H}_{s,d} (\vb*{d}_2)$ is still open: we did not find examples in this space where one of the methods samples uniformly and the other samples with bias. 

\begin{table}[tbp]
\centering
\begin{threeparttable}
    \centering
    \caption{Hypercurveball for undirected hypergraphs.}
    \begin{tabular}{l|c|c||c|c}
        Configuration & Degenerate & Multi- & Always & Always uniform\\ 
        model & hyperedges & hyperedges & uniform? & for graphs? \\ \hline
        $\CM{}{d}{1}$ & No & No & ?\tnote{1} & Yes \cite{carstens2016}\\
        $\CM{d}{d}{1}$ & Yes & No & No & No\\
        $\CM{m}{d}{1}$ & No & Yes & Yes & Yes\\
        $\CM{d,m}{d}{1}$ & Yes & Yes & Yes & Yes  
    \end{tabular}
     \begin{tablenotes}
     \item[1] No counterexamples exist for $|V|,|E| \leq 4$.
  \end{tablenotes}
\label{tab:results_undir}
\end{threeparttable}
\end{table}

\begin{table}[tbp]
\centering
\begin{threeparttable}
    \centering
    \caption{Hypercurveball for directed hypergraphs.}
    \begin{tabular}{l|c|c|c||c|c}
        Configuration & Self-loops & Degenerate & Multiple & Always & Always uniform\\ 
        model & & hyperarcs & hyperarcs & uniform? & for digraphs?\\ \hline
        $\CM{}{d}{2}$ & No & No & No & ?\tnote{1} & Yes\\
        $\CM{d}{d}{2}$ & No & Yes & No & No & -\\
        $\CM{m}{d}{2}$ & No & No & Yes & ?\tnote{2}& Yes\\
        $\CM{d,m}{d}{2}$ & No & Yes & Yes & No & - \\
        $\CM{s}{d}{2}$ & Yes & No & No & ?\tnote{3} & Yes \cite{carstens2015}\\
        $\CM{s,d}{d}{2}$ & Yes & Yes & No & No & - \\
        $\CM{s,m}{d}{2}$ & Yes & No & Yes & Yes & Yes\\
        $\CM{s,d,m}{d}{2}$ & Yes & Yes & Yes & Yes & -       
    \end{tabular}
     \begin{tablenotes}
    \item[1] No counterexamples exist for $|V|\leq 4,|E| \leq 3$, nor for $|V|\leq 3,|E|=4$.
  \item[2] No counterexamples exist for $|V|\leq 4, |E| \leq 3$, nor for $|V|=5,|E|=2$.
  \item[3] No counterexamples exist for $|V| \leq 4, |E| \leq 2, \forall v \in V, e \in E:  m_{e^{\tail}}(v),m_{e^{\head}}(v) \leq 2 $.
  \end{tablenotes}
\label{tab:results_dir}
\end{threeparttable}
\end{table}

\subsection{Comparison to the hyperedge-shuffle method}
\label{section:comparison}

In the remainder of this work, we compare the Hypercurveball method to the hyperedge-shuffle method \cite{chodrow2019,kraakman2024}. This subsection focuses on comparing the theoretical results of the two methods, while Section \ref{section:mixing} covers experimental results.

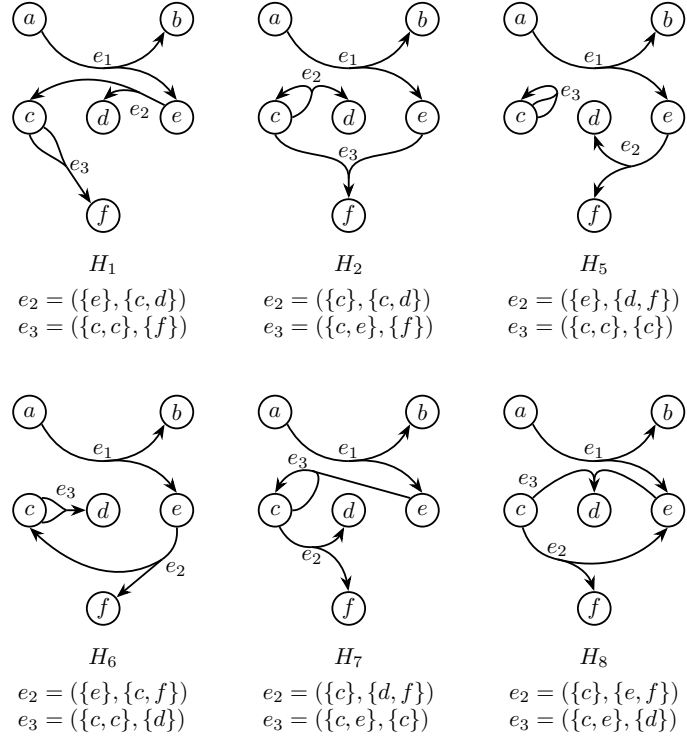
\begin{figure}[tbp]
    \centering
    \newcommand{\offsetx}{5*\tikzscale}
\newcommand{\offsety}{-8*\tikzscale}

\resizebox{0.6\textwidth}{!}{%
\begin{tikzpicture}[-Stealth, line width=1.3pt,auto,
                    thick,main node/.style={circle,draw, minimum size=\tikznodesize}, inner sep=1pt]

  \node[main node] (a1) at (0,4*\tikzscale) {$a$};
  \node[main node] (b1) at (3*\tikzscale,4*\tikzscale) {$b$};
  \node[main node] (c1) at (0,2*\tikzscale) {$c$};
  \node[main node] (d1) at (1.5*\tikzscale,2*\tikzscale) {$d$};
  \node[main node] (e1) at (3*\tikzscale,2*\tikzscale) {$e$};
  \node[main node] (f1) at (1.5*\tikzscale,0) {$f$};

  \node (x) at (1.5*\tikzscale, -\tikzscale) {$H_1$};
  \node[align=left] at (1.5*\tikzscale, -2*\tikzscale) {$e_2= (\{e\},\{c,d\}) $\\$e_3 = (\{c,c\},\{f\})$};

\coordinate (abe1) at (1.5*\tikzscale,3*\tikzscale); 

\draw[-, bend left] (abe1) to (a1.315);
\draw[bend right] (abe1) to (b1.225);
\draw[bend left] (abe1) to (e1.90);

\node[above] at (abe1) {$e_1$};

\coordinate (ecd1) at (2.25*\tikzscale,2.5*\tikzscale); 

\draw[-] (ecd1) to (e1.135);
\draw[bend right] (ecd1) to (c1.90);
\draw[bend right] (ecd1) to (d1.90);

\node[below, yshift=-5*\tikzscale pt] at (ecd1) {$e_2$};

\coordinate (ccf1) at (0.75*\tikzscale,\tikzscale); 

\draw[-,bend right, out=0] (ccf1) to (c1.325);
\draw[-,bend left, out=0] (ccf1) to (c1.270);
\draw (ccf1) to (f1.135);

\node[right] at (ccf1) {$e_3$};

 \node[main node] (a2) at (\offsetx,4*\tikzscale) {$a$};
  \node[main node] (b2) at (\offsetx+3*\tikzscale,4*\tikzscale) {$b$};
  \node[main node] (c2) at (\offsetx,2*\tikzscale) {$c$};
  \node[main node] (d2) at (\offsetx+1.5*\tikzscale,2*\tikzscale) {$d$};
  \node[main node] (e2) at (\offsetx+3*\tikzscale,2*\tikzscale) {$e$};
  \node[main node] (f2) at (\offsetx+1.5*\tikzscale,0) {$f$};

  \node (x) at (\offsetx + 1.5*\tikzscale, -\tikzscale) {$H_2$};
\node[align=left] at (\offsetx + 1.5*\tikzscale, -2*\tikzscale) {$e_2= (\{c\},\{c,d\}) $\\$e_3 = (\{c,e\},\{f\})$};

\coordinate (abe2) at (\offsetx+1.5*\tikzscale,3*\tikzscale); 

\draw[-, bend left] (abe2) to (a2.315);
\draw[bend right] (abe2) to (b2.225);
\draw[bend left] (abe2) to (e2.90);

\node[above] at (abe2) {$e_1$};

\coordinate (ccd2) at (\offsetx+ 0.75*\tikzscale,2.5*\tikzscale); 

\draw[-, bend left] (ccd2) to (c2.0);
\draw[bend right, out=270] (ccd2) to (c2.90);
\draw[bend left, out=90] (ccd2) to (d2.90);

\node[above, yshift = 4.25 * \tikzscale pt] at (ccd2) {$e_2$};

\coordinate (cef2) at (\offsetx + 1.5*\tikzscale,0.75*\tikzscale); 

\draw[-,bend left, out=300] (cef2) to (c2.270);
\draw[-,bend right, out=60] (cef2) to (e2.270);
\draw (cef2) to (f2.90);

\node[above, yshift = 9*\tikzscale pt] at (cef2) {$e_3$};

  \node[main node] (a3) at (2*\offsetx,4*\tikzscale) {$a$};
  \node[main node] (b3) at (2*\offsetx+3*\tikzscale,4*\tikzscale) {$b$};
  \node[main node] (c3) at (2*\offsetx,2*\tikzscale) {$c$};
  \node[main node] (d3) at (2*\offsetx+1.5*\tikzscale,2*\tikzscale) {$d$};
  \node[main node] (e3) at (2*\offsetx+3*\tikzscale,2*\tikzscale) {$e$};
  \node[main node] (f3) at (2*\offsetx+1.5*\tikzscale,0) {$f$};

  \node (x) at (2*\offsetx+1.5*\tikzscale, -\tikzscale) {$H_5$};
  \node[align=left] at (2*\offsetx + 1.5*\tikzscale, -2*\tikzscale) {$e_2= (\{e\},\{d,f\}) $\\$e_3 = (\{c,c\},\{c\})$};

\coordinate (abe3) at (2*\offsetx+1.5*\tikzscale,3*\tikzscale); 

\draw[-, bend left] (abe3) to (a3.315);
\draw[bend right] (abe3) to (b3.225);
\draw[bend left] (abe3) to (e3.90);

\node[above] at (abe3) {$e_1$};

\coordinate (ccc3) at (2*\offsetx+0.75*\tikzscale, 2.5*\tikzscale); 

\draw[-, bend left, in=135] (ccc3) to (c3.0);
\draw[-, bend left, in=200] (ccc3) to (c3.45);
\draw[bend right, out=270] (ccc3) to (c3.90);

\node[right] at (ccc3) {$e_3$};

\coordinate (ecf3) at (2*\offsetx+2.25*\tikzscale,\tikzscale); 

\draw[-,bend right] (ecf3) to (e3.270);
\draw[bend left] (ecf3) to (d3.270);
\draw[bend right] (ecf3) to (f3.90);

\node[above, yshift = 5*\tikzscale pt] at (ecf3) {$e_2$};

 \node[main node] (a4) at (0,\offsety + 4*\tikzscale) {$a$};
  \node[main node] (b4) at (3*\tikzscale,\offsety + 4*\tikzscale) {$b$};
  \node[main node] (c4) at (0,\offsety + 2*\tikzscale) {$c$};
  \node[main node] (d4) at (1.5*\tikzscale,\offsety+ 2*\tikzscale) {$d$};
  \node[main node] (e4) at (3*\tikzscale,\offsety + 2*\tikzscale) {$e$};
  \node[main node] (f4) at (1.5*\tikzscale,\offsety) {$f$};

  \node (x) at (1.5*\tikzscale,\offsety -\tikzscale) {$H_6$};
  \node[align=left] at (1.5*\tikzscale, \offsety -2*\tikzscale) {$e_2= (\{e\},\{c,f\}) $\\$e_3 = (\{c,c\},\{d\})$};

\coordinate (abe4) at (1.5*\tikzscale,\offsety + 3*\tikzscale); 

\draw[-, bend left] (abe4) to (a4.315);
\draw[bend right] (abe4) to (b4.225);
\draw[bend left] (abe4) to (e4.90);

\node[above] at (abe4) {$e_1$};

\coordinate (ccd4) at (0.75*\tikzscale,\offsety + 2*\tikzscale); 

\draw[-, bend right, out=0] (ccd4) to (c4.45);
\draw[-, bend left, out=0] (ccd4) to (c4.315);
\draw (ccd4) to (d4.180);

\node[above, yshift = 5*\tikzscale pt] at (ccd4) {$e_3$};

\coordinate (ecf4) at (2.7*\tikzscale,\offsety + \tikzscale); 

\draw[-, bend right] (ecf4) to (e4.270);
\draw[bend left, out=45] (ecf4) to (c4.270);
\draw[] (ecf4) to (f4);

\node[below right] at (ecf4) {$e_2$};

\node[main node] (a5) at (\offsetx,\offsety + 4*\tikzscale) {$a$};
  \node[main node] (b5) at (\offsetx+3*\tikzscale,\offsety + 4*\tikzscale) {$b$};
  \node[main node] (c5) at (\offsetx,\offsety + 2*\tikzscale) {$c$};
  \node[main node] (d5) at (\offsetx+1.5*\tikzscale,\offsety + 2*\tikzscale) {$d$};
  \node[main node] (e5) at (\offsetx+3*\tikzscale,\offsety + 2*\tikzscale) {$e$};
  \node[main node] (f5) at (\offsetx+1.5*\tikzscale,\offsety) {$f$};

  \node (x) at (\offsetx + 1.5*\tikzscale,\offsety -\tikzscale) {$H_7$};
    \node[align=left] at (\offsetx + 1.5*\tikzscale, \offsety -2*\tikzscale) {$e_2= (\{c\},\{d,f\}) $\\$e_3 = (\{c,e\},\{c\})$};

\coordinate (abe5) at (\offsetx+1.5*\tikzscale,\offsety + 3*\tikzscale); 

\draw[-, bend left] (abe5) to (a5.315);
\draw[bend right] (abe5) to (b5.225);
\draw[bend left] (abe5) to (e5.90);

\node[above] at (abe5) {$e_1$};

\coordinate (cec5) at (\offsetx+ 0.75*\tikzscale,\offsety + 2.8*\tikzscale); 

\draw[-, bend left, out=115, in=110] (cec5) to (c5.0);
\draw[bend right, out=315] (cec5) to (c5.90);
\draw[-] (cec5) to (e5.135);

\node[above left] at (cec5) {$e_3$};

\coordinate (cdf5) at (\offsetx + 0.75*\tikzscale,\offsety + 1.25*\tikzscale); 

\draw[-, bend left, out=35] (cdf5) to (c5);
\draw[bend right, out=325] (cdf5) to (d5);
\draw[bend left, out=45] (cdf5) to (f5);

\node[below] at (cdf5) {$e_2$};

\node[main node] (a6) at (2*\offsetx,\offsety + 4*\tikzscale) {$a$};
  \node[main node] (b6) at (2*\offsetx+3*\tikzscale,\offsety + 4*\tikzscale) {$b$};
  \node[main node] (c6) at (2*\offsetx,\offsety + 2*\tikzscale) {$c$};
  \node[main node] (d6) at (2*\offsetx+1.5*\tikzscale,\offsety + 2*\tikzscale) {$d$};
  \node[main node] (e6) at (2*\offsetx+3*\tikzscale,\offsety + 2*\tikzscale) {$e$};
  \node[main node] (f6) at (2*\offsetx+1.5*\tikzscale,\offsety) {$f$};

  \node (x) at (2*\offsetx + 1.5*\tikzscale, \offsety-\tikzscale) {$H_8$};
\node[align=left] at (2*\offsetx + 1.5*\tikzscale, \offsety -2*\tikzscale) {$e_2= (\{c\},\{e,f\}) $\\$e_3 = (\{c,e\},\{d\})$};

\coordinate (abe6) at (2*\offsetx+1.5*\tikzscale,\offsety +  3*\tikzscale); 

\draw[-, bend left] (abe6) to (a6.315);
\draw[bend right] (abe6) to (b6.225);
\draw[bend left] (abe6) to (e6.90);

\node[above] at (abe6) {$e_1$};

\coordinate (ced6) at (2*\offsetx+ 1.5*\tikzscale,\offsety+2.6*\tikzscale); 

\draw[-, out=90, in=45] (ced6) to (c6.45);
\draw[-, out=90] (ced6) to (e6);
\draw[] (ced6) to (d6);

\node[left, xshift=-30*\tikzscale pt] at (ced6) {$e_3$};

\coordinate (cef6) at (2*\offsetx + 0.75*\tikzscale,\offsety+\tikzscale); 

\draw[-,bend left] (cef6) to (c6);
\draw[bend right] (cef6) to (e6.270);
\draw[bend left] (cef6) to (f6.90);

\node[above] at (cef6) {$e_2$};

\end{tikzpicture}
}
    \caption{Performing a hyperedge-shuffle on the hyperedges $e_2$ and $e_3$ in $H_1$ (same as in Figure \ref{fig:example_trade}) can result in the hypergraph $H_1,H_2,H_5,H_6,H_7$ or $H_8$.}
    \label{fig:example_shuffle}
\end{figure}

Firstly, we introduce the hyperedge-shuffle method \cite{chodrow2019,kraakman2024}. At each step of this method for hypergraphs, two hyperedges $e,f \in E$ are chosen uniformly at random. Then, the elements of these edges are shuffled. More precisely, $e \uplus f$ is randomly partitioned into multisets $e', f'$ with cardinalities $|e|$ and $|f|$, respectively. For directed hypergraphs, the sets $e^{\tail} \uplus f^{\tail}$ and $e^{\head} \uplus f^{\head}$ are randomly partitioned, similar to how the sets $I^{\tail}_v, I^{\tail}_w, I^{\head}_v$ and $I^{\head}_w$ are used in the hypercurveball algorithm. To illustrate, consider the directed hypergraph $H_1$ in Figure \ref{fig:example_trade}. Performing a hyperedge-shuffle on the hyperedges $(\{c,c\},\{f\})$ and $(\{e\},\{c,d\})$ in this hypergraph can result in the hypergraph $H_1, H_2$ or the hypergraph $H_5,H_6,H_7$ or $H_8$ as in Figure \ref{fig:example_shuffle}.

In Section \ref{section:curveball}, we explained the intuition of why the Curveball method mixes faster on graphs than the edge-swapping method: only four nodes can be incident to two edges, but many edges can be incident to two nodes. Therefore, the Curveball method can reconfigure more elements in one time step. This argument does not hold for hypergraphs since hyperedges can contain many nodes. Interestingly, we show in Section \ref{section:mixing} that the Hypercurveball method is indeed not always faster than the hyperedge-swapping method; for some initial hypergraphs, the Hypercurveball method is faster, while for others, the hyperedge-swapping method prevails.

\begin{figure}[tbp]
    \centering
    \begin{subfigure}{0.4\textwidth}
    \centering
    \begin{tikzpicture}[-Stealth,line width=1.3pt,auto,
                    thick,main node/.style={circle,draw, minimum size=\tikznodesize}, inner sep=1pt,every loop/.style={}]

  \node[main node] (u) at (0,0) {$u$};
  \node[main node] (v) at (3*\tikzscale,0) {$v$};

\coordinate (vvuu) at (1.5*\tikzscale,0);

\draw[-,bend right, in = 170] (v.135) to (vvuu);
\draw[-,bend left, in=190] (v.225) to (vvuu);
\draw[bend right,out=10] (vvuu) to (u.45);
\draw[bend left, out=350] (vvuu) to (u.315);

\draw[bend left] (u.90) to (v.90);

\coordinate (uuuv) at (0,-\tikzscale);

\draw[-,bend right, in=255] (u.250) to (uuuv);
\draw[-,bend right, out=300, in=255] (u.180) to (uuuv);
\draw[bend right, out=280] (uuuv) to (u.290);
\draw[bend right, out=350] (uuuv) to (v.270);

\end{tikzpicture}
        \caption{$H_1$}
    \end{subfigure}
    \begin{subfigure}{0.4\textwidth}
    \centering
    \begin{tikzpicture}[-Stealth,line width=1.3pt,auto,
                    thick,main node/.style={circle,draw, minimum size=\tikznodesize}, inner sep=1pt,every loop/.style={}]

  \node[main node] (u) at (0,0) {$u$};
  \node[main node] (v) at (3*\tikzscale,0) {$v$};

\coordinate (uuvv) at (1.5*\tikzscale,0);

\draw[-,bend right, in = 170] (u.315) to (uuvv);
\draw[-,bend left, in=190] (u.45) to (uuvv);
\draw[bend right,out=10] (uuvv) to (v.225);
\draw[bend left, out=350] (uuvv) to (v.135);

\draw[bend right] (v.90) to (u.90);

\coordinate (uvuu) at (0,-\tikzscale);

\draw[-, bend right, out=280] (uuuv) to (u.290);
\draw[-, bend right, out=350] (uuuv) to (v.270);

\draw[bend left, out=75] (uuuv) to (u.250);
\draw[bend left, in=120, out=75] (uuuv) to (u.180);

\end{tikzpicture}
        \caption{$H_2$}
    \end{subfigure}
    \caption{Two directed hypergraphs $H_1,H_2 \in (\mathcal{H}_{d,m}(\vb*{d}^{**}) \cap \mathcal{H}_{d}(\vb*{d}^{**}))$ which are not connected using hyperedge-shuffles, but which are connected using hypertrades in either space $\mathcal{H}_{d,m}(\vb*{d}^{**})$ or $\mathcal{H}_d(\vb*{d}^{**})$.}
    \label{fig:noneq_fig2}
\end{figure}
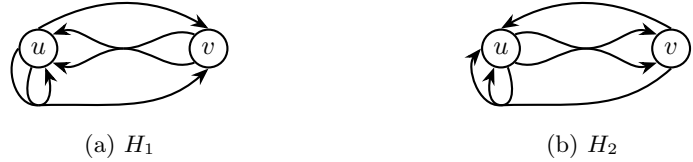

Next, we address the symmetry of the two methods. The hyperedge-shuffle procedure and the hypertrade are quite similar when considering the incidence matrix of a hypergraph (with a row for each node and a column for each hyperedge). With a hyperedge-shuffle, the entries in two columns are randomly exchanged, and with a hypertrade, the entries in two rows are randomly exchanged. Still, the methods are not equivalent in terms of uniform sampling (Theorem \ref{thm:methods_not_equivalent}) or in terms of the three special hyperedge types $s$, $d$ and $m$, because the hyperedge types refer specifically to the columns of the incidence matrix. Therefore, different restrictions are imposed on the hyperedge-shuffle and the hypertrade. For instance, a hyperedge-shuffle can create at most two degenerate hyperedges, namely the two hyperedges that were shuffled, but a hypertrade can create many degenerate hyperedges. In addition, consider the two hypergraphs in Figure \ref{fig:biased_dir_H_dm}. With one hyperedge-shuffle, $H_1$ can be transformed into $H_2$, but it is impossible to transform $H_1$ into $H_2$ using hypertrades without creating a self-loop, which is not allowed. On the contrary, consider the two hypergraphs in Figure \ref{fig:noneq_fig2}. With one hypertrade, $H_1$ can be transformed into $H_2$, but it is impossible to transform $H_1$ into $H_2$ using hypertrades without creating a self-loop, which is not allowed. For more details of these examples, see the proof of Theorem \ref{thm:methods_not_equivalent} in Section \ref{section:pf_thm4}.

\begin{figure}[tbp]
    \centering
    \begin{subfigure}{0.4\textwidth}
    \centering
    \begin{tikzpicture}[-Stealth,line width=1.3pt,auto,
                    thick,main node/.style={circle,draw, minimum size=\tikznodesize}, inner sep=1pt,every loop/.style={}]

  \node[main node] (v) at (0,0) {$v$};
  \node[main node] (w) at (2*\tikzscale,0) {$w$};
  \node[main node] (u) at (\tikzscale,1.73*\tikzscale) {$u$};

\coordinate (uvuu) at (0,1.5*\tikzscale);

\draw[-,bend left] (u.225) to (uvuu);
\draw[-,bend right] (v.90) to (uvuu);
\draw[bend left, out=110, in=160] (uvuu) to (u.180);
\draw[bend left, out=90, in=140] (uvuu) to (u.135);

\draw[] (u) to (w);

\end{tikzpicture}
        \caption{$H_1$}
    \end{subfigure}
    \begin{subfigure}{0.4\textwidth}
    \centering
    \begin{tikzpicture}[-Stealth,line width=1.3pt,auto,
                    thick,main node/.style={circle,draw, minimum size=\tikznodesize}, inner sep=1pt,every loop/.style={}]

  \node[main node] (v) at (0,0) {$v$};
  \node[main node] (w) at (2*\tikzscale,0) {$w$};
  \node[main node] (u) at (\tikzscale,1.73*\tikzscale) {$u$};

\coordinate (uuuw) at (2*\tikzscale,1.5*\tikzscale);

\draw[-,bend left, out=45, in=90] (u.45) to (uuuw);
\draw[-,bend left, in=80] (u.0) to (uuuw);
\draw[bend left] (uuuw) to (u.300);
\draw[bend right] (uuuw) to (w);

\draw[] (v) to (u);

\end{tikzpicture}
        \caption{$H_2$}
    \end{subfigure}
    \caption{Two directed hypergraphs $H_1,H_2 \in (\mathcal{H}_{d,m}(\vb*{d}^*) \cap \mathcal{H}_{d}(\vb*{d}^*))$ which are connected using hyperedge-shuffles, and not connected using hypertrades in either space $\mathcal{H}_{d,m}(\vb*{d}^*)$ or $\mathcal{H}_d(\vb*{d}^*)$.}
    \label{fig:biased_dir_H_dm}
\end{figure}
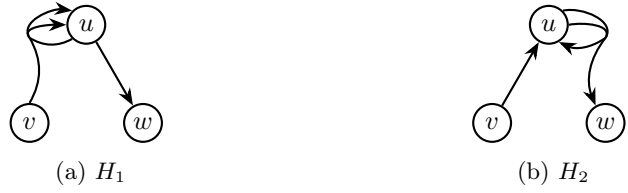

\newpage
\section{Mixing time}
\label{section:mixing}
When applying the Hypercurveball algorithm, the number of trades to perform, $N$, must be selected (see Algorithm \ref{alg:curveball_sdm}). This number should reflect the mixing time of the method, i.e., the number of steps of the algorithm needed to obtain a uniformly random hypergraph. This section experimentally compares the mixing time of the Hypercurveball algorithm to that of the hyperedge-shuffle algorithm. 

To experimentally demonstrate the mixing behavior of both sampling methods on hypergraphs, we use the \textit{perturbation degree} of the hypergraph after each step of the algorithm relative to the initial hypergraph. The perturbation degree measures the dissimilarity between two hypergraphs, based on the definition in \cite{strona2014}. We use the point at which the perturbation degree stabilizes as an estimate of the mixing time. Specifically, we estimate the mixing time as the first step at which the perturbation degree is within a 2\% range of $L$, where $L$ is the limiting value of the perturbation degree.

The perturbation degree of a hypergraph $H$ with respect to the initial hypergraph $G$, $p_G(H)$, is computed by first mapping each hyperedge $e$ in $E(H)$ to a hyperedge $\Tilde{e}$ in $G$ with the same degree. The perturbation degree is defined as follows.
\\
\begin{definition}[Perturbation degree]
\label{def:perturbation}
    Let $G,H \in \mathcal{H}_{x} (\vb*{d})$. The perturbation degree of $H$ w.r.t. $G$ is   
    \begin{align*}
        p_G(H) = \frac{1}{2\sum_{v \in V} d_v} \sum_{v \in V} \sum_{e \in E(H)} |m_e(v) - m_{\Tilde{e}}(v)|.
    \end{align*}
    
\end{definition}
For directed hypergraphs, the perturbation degree is calculated by summing the perturbation degree of the tails and the heads. Appendix \ref{app:perturbation} contains the definition of the perturbation degree for directed hypergraphs and implementation details. 

We analyze five datasets, of which two are undirected hypergraphs and three are directed hypergraphs, see Table~\ref{tab:datasets}:
\begin{itemize}
    \item \texttt{unicodelang} \cite{kunegis2013}: language spoken by country. Nodes: countries. Hyperedges: languages. Each hyperedge contains the set of countries that speak the language associated to the hyperedge.
    \item \texttt{board\_directors} \cite{seierstad2011}: Norwegian boards of directors. Nodes: people. Hyperedges: boards. Each hyperedge contains the set of people that are in the board together.
    \item \texttt{email\_eu} \cite{leskovec2007}: email network. Nodes: people. Directed hyperedges: emails. Each directed hyperedge contains the email sender (tail) and the set of receivers (head).
    \item \texttt{thiol\_oscillator} \cite{semenov2016}: chemical reaction network. Nodes: chemical species. Directed hyperedges: reactions. Each directed hyperedge contains the multiset of reactants (tail) and products (head) of the reaction.
    \item \texttt{iAF1260b} \cite{yadati2020}: metabolic reactions. Nodes: chemical species. Directed hyperedges: reactions. Each directed hyperedge contains the multiset of reactants (tail) and products (head) of the reaction.
\end{itemize}

We then investigate the influence of the following statistics of the datasets on the mixing time of the Hypercurveball method:
\begin{itemize}
    \item The number of nodes. If the number of nodes is small, then each hypertrade can modify the hypergraph by a larger proportion than when the number of nodes is larger.
    \item The average node degree. If the average node degree is large, then an average hypertrade can change more hyperedges in one time step.
    \item The median node degree: Similarly to the average node degree, a typical hypertrade can change more hyperedges in one time step when the median is larger.
    \item $\mathds{E}_V[\min \{d_1,d_2\}]$. This is the expected minimum degree of two randomly chosen nodes. Each hypertrade can change the incidence sets of the chosen nodes $u$ and $v$ in as many places as the minimum of the two node degrees (Algorithm \ref{alg:trade_d}). Therefore, the bigger this minimum, the more a hypertrade can modify the hypergraph. Appendix \ref{app:minimum_of_two} shows how to compute this statistic.
\end{itemize}
The same statistics for the hyperedges are investigated as well, which should influence the mixing time of the hyperedge-shuffle method in a similar way. For each dataset, its statistics are summarized in Table \ref{tab:datasets}. 

\begin{table}[tbp]
\centering
\begin{threeparttable}
    \centering
    \caption{Data sets used for experiments on mixing time. Lines with tuple data refer to data modeled as directed hypergraphs, where the first entry in the tuple shows the statistics for the tails and the second entry for the heads.}
    \begin{tabular}{|c c|c|c|c|c|}
    \hline
       Data set & & \# & Average degree & Median degree & $\mathds{E}[\min\{d_1,d_2\}]$\tnote{1}\\ \hline \hline
       \texttt{unicodelang} & nodes & 254 & 4.94 & 3.00 & 2.34 \\
       \cite{kunegis2013} & hyperedges & 614 & 2.04 & 1.00 & 1.08 \\ \hline
       \texttt{board\_directors} & nodes & 1013 & 1.12 & 1.00 & 1.01 \\
        \cite{seierstad2011} & hyperedges & 204 & 5.54 & 5.00 & 4.32 \\ \hline
    \texttt{email\_eu} & nodes & 986 & (36.28, 85.90) & (16.00, 46.00) & (12.00, 33.79) \\
    \cite{leskovec2007} & hyperedges & 35774 & (1.00, 2.37) & (1.00, 1.00) & (1.00, 1.30)\\ \hline
    \texttt{thiol\_oscillator} & nodes & 18 & (4.94, 4.72) & (4.00, 2.50) & (2.75, 2.21) \\
    \cite{semenov2016}& hyperedges & 55 & (1.62, 1.55) & (2.00, 2.00) & (1.38, 1.27) \\ \hline
    \texttt{iAF1260b} & nodes & 1668 & (2.53, 2.83) & (1.00, 1.00) & (0.97, 0.76) \\ 
    \cite{yadati2020} & hyperedges & 2084 & (2.03, 2.27) & (2.00, 2.00) & (1.62, 1.77) \\ \hline
    \end{tabular}
     \begin{tablenotes}
    \item[1] The expected minimum degree of two randomly chosen nodes/hyperedges.
  \end{tablenotes}
\label{tab:datasets}
\end{threeparttable}
\end{table}

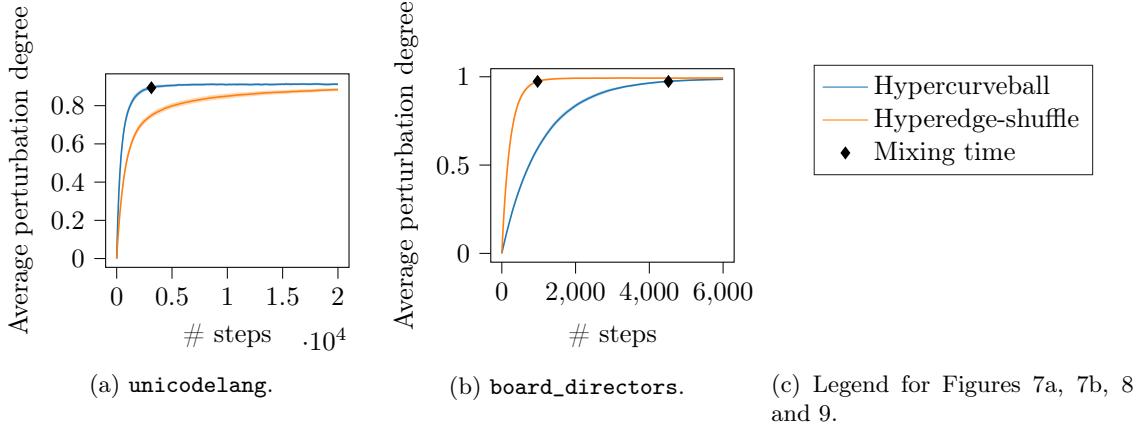
\begin{figure}[tbp]
\centering
\begin{subfigure}[t]{0.32\textwidth}
    \centering
    \subfile{Tikz/MP_unicodelang_reduced.tex}
    \caption{\texttt{unicodelang}.}
    \label{fig:experiment_undir_1}
\end{subfigure} \hfill
\begin{subfigure}[t]{0.32\textwidth}
    \centering
    \subfile{Tikz/MP_board_directors_reduced.tex}
    \vspace{-0.4cm}
    \caption{\texttt{board\_directors}.}
    \label{fig:experiment_undir_2}
\end{subfigure} \hfill
\begin{subfigure}[t]{0.32\textwidth}
    \centering
    \begin{tikzpicture} 

\definecolor{darkgray176}{RGB}{176,176,176}
\definecolor{darkorange25512714}{RGB}{255,127,14}
\definecolor{lightgray204}{RGB}{204,204,204}
\definecolor{steelblue31119180}{RGB}{31,119,180}

    \begin{axis}[%
    hide axis,
    xmin=10,
    xmax=50,
    ymin=0,
    ymax=0.4,
    legend style={draw=white!15!black,legend cell align=left}
    ]
    \addlegendimage{steelblue31119180}
    \addlegendentry{Hypercurveball};
    \addlegendimage{darkorange25512714}
    \addlegendentry{Hyperedge-shuffle};
    \addlegendimage{only marks, semithick, black, mark = diamond*}
    \addlegendentry{Mixing time}
    \end{axis}

\node[] () [below = 1.5em] at (3,2.5) {}; %

\end{tikzpicture}
    \caption{Legend for Figures \ref{fig:experiment_undir_1}, \ref{fig:experiment_undir_2}, \ref{fig:experiment_dir} and \ref{fig:experiment_artificial}.}
    \label{fig:mixing_plot_legend}
\end{subfigure} 
\caption{Experimental results for the mixing time of two undirected hypergraph data sets. Plots are the average of 100 runs of the two algorithms on $\mathcal{H}_{d,m}$. The shaded area is one standard deviation. For the hyperedge-shuffle method on the \texttt{unicodelang} dataset, the mixing time lies outside of our experimental range. The legend for Figures \ref{fig:experiment_undir_1} and \ref{fig:experiment_undir_2} is shown in Figure \ref{fig:mixing_plot_legend}.}
\label{fig:experiment_undir}
\end{figure}

\begin{figure}[tbp]
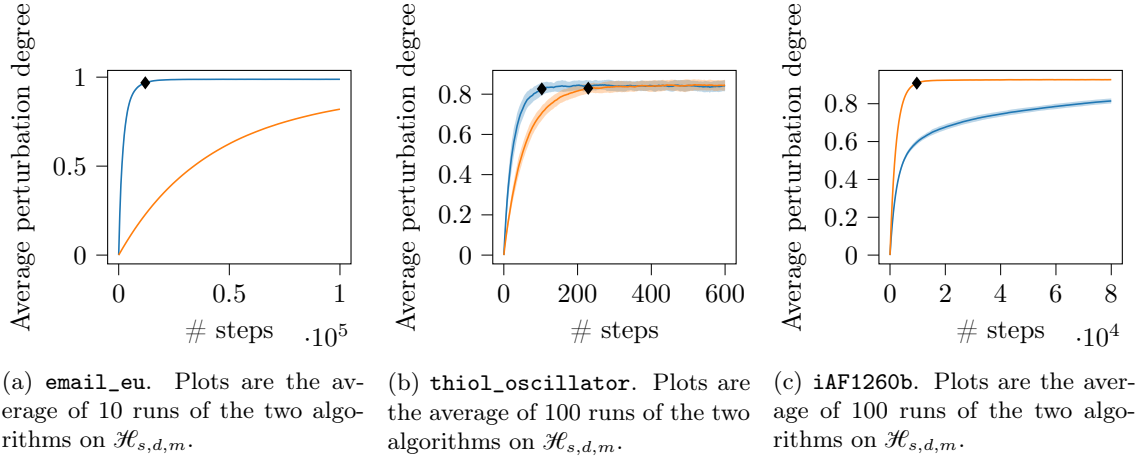

\centering
\begin{subfigure}[t]{0.32\textwidth}
    \centering
    \subfile{Tikz/MP_email_eu_reduced.tex}
    \caption{\texttt{email\_eu}. Plots are the average of 10 runs of the two algorithms on $\mathcal{H}_{s,d,m}$.}
\end{subfigure} \hfill
\begin{subfigure}[t]{0.32\textwidth}
    \centering
    \subfile{Tikz/MP_thiol_oscillator_reduced.tex}
    \vspace{-0.4cm}
    \caption{\texttt{thiol\_oscillator}. Plots are the average of 100 runs of the two algorithms on $\mathcal{H}_{s,d,m}$.}
\end{subfigure} \hfill
\begin{subfigure}[t]{0.32\textwidth}
    \centering
    \subfile{Tikz/MP_iAF1260b_reduced.tex}
    \caption{\texttt{iAF1260b}. Plots are the average of 100 runs of the two algorithms on $\mathcal{H}_{s,d,m}$.}
\end{subfigure}
\caption{Experimental results for the mixing time of three directed hypergraph data sets. The shaded area is one standard deviation. For the hyperedge-shuffle method on the \texttt{email\_eu} dataset and the Hypercurveball method on the \texttt{iAF1260b} dataset, the mixing time lies outside of our experimental range. The legend for this figure is shown in Figure \ref{fig:mixing_plot_legend}.}
\label{fig:experiment_dir}
\end{figure}

We will now show how, given a data set, these statistics for the nodes and for the hyperedges can indicate whether the Hypercurveball method or the hyperedge-shuffle method mixes faster.
The experimental mixing times of the two undirected and directed hypergraph data sets are shown in Figures \ref{fig:experiment_undir} and \ref{fig:experiment_dir} respectively.

The results reveal that some datasets favor the Hypercurveball method for faster mixing, while others favor the hyperedge-shuffle method. Comparing Table \ref{tab:datasets} with the experimental outcomes, it is evident that merely comparing the number of nodes and hyperedges does not determine which method is the fastest. For instance, in the \texttt{unicodelang} dataset, which has more hyperedges than nodes, the Hypercurveball method mixes faster, whereas in the \texttt{iAF1260b} dataset with similar characteristics, the hyperedge-shuffle method is faster. The same observation holds for the average degrees.

\begin{figure}[tbp]
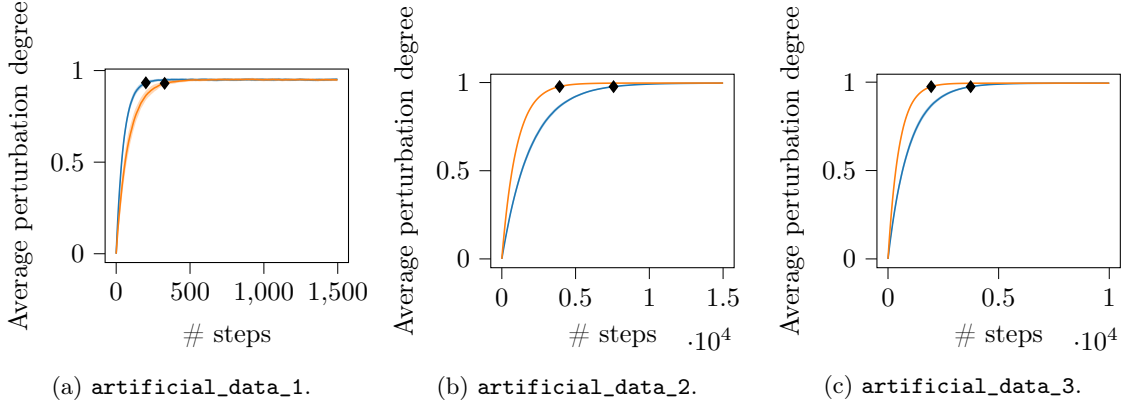

\centering
\begin{subfigure}[t]{0.32\textwidth}
    \centering
    \subfile{Tikz/MP_artificial_1_reduced.tex}
    \vspace{-0.4cm}
    \caption{\texttt{artificial\_data\_1}.}
\end{subfigure} \hfill
\begin{subfigure}[t]{0.32\textwidth}
    \centering
    \subfile{Tikz/MP_artificial_2_reduced.tex}
    \vspace{-0.4cm}
    \caption{\texttt{artificial\_data\_2}.}
\end{subfigure}
\hfill
\begin{subfigure}[t]{0.32\textwidth}
    \centering
    \subfile{Tikz/MP_artificial_3_reduced.tex}
    \caption{\texttt{artificial\_data\_3}.}
\end{subfigure}
\caption{Experimental results for the mixing time of three artificial data sets. Plots are the average of 100 runs of the two algorithms on $\mathcal{H}_{s,d,m}$. The shaded area is one standard deviation. The legend for this figure is shown in Figure \ref{fig:mixing_plot_legend}.}
\label{fig:experiment_artificial}
\end{figure}

On the other hand, the results indicate that the Hypercurveball method always mixes faster than the hyperedge-shuffle method when the median of the node degrees exceeds the median of the hyperedge degrees. This finding also applies to the expected minimum degree of two nodes/hyperedges. To delve deeper into these statistics, we analyze three artificial datasets: one where the median node degree is smaller than the median hyperedge degree, and the expected minimum degree of two nodes is larger than that of two hyperedges, and two where the reverse is true. See Table \ref{tab:artificial_dataset} and Appendix \ref{app:artificial_dataset} for the statistics and the details of these datasets. The experimental results for the mixing time are shown in Figure \ref{fig:experiment_artificial}.

\begin{table}[tbp]
\centering
\begin{threeparttable}
    \centering
    \caption{Artificial data sets used for experiments on mixing time.}
    \begin{tabular}{|c c|c|c|c|c|}
    \hline
       Data set & & \# & Average degree & Median degree & $\mathds{E}[\min\{d_1,d_2\}]$\tnote{1}\\ \hline \hline
       \texttt{artificial\_data\_1} & nodes & 51 & 16.27 & 16.00 & 16.00 \\
        & hyperedges & 51 & 16.27 & 30.00 & 9.14 \\ \hline
        \texttt{artificial\_data\_2} & nodes & 1001 & 50.05 & 99.00 & 25.52 \\
        & hyperedges & 1001 & 50.05 & 50.00 & 50.00 \\
        \hline
        \texttt{artificial\_data\_3} & nodes & 501 & 25.05 & 49.00 & 13.02 \\
        & hyperedges & 501& 25.05 & 25.00 & 25.00 \\ \hline
    \end{tabular}
     \begin{tablenotes}
    \item[1] The expected minimum degree of two nodes/hyperedges.
  \end{tablenotes}
\label{tab:artificial_dataset}
\end{threeparttable}
\end{table}

\begin{figure}[tbp]
    \centering
    \subfile{Tikz/MP_all_data_reduced}
    \begin{tikzpicture}
    \begin{axis}[
    width = 0.3*\linewidth,
        hide axis,
        scale only axis,
        height=6.2cm,
        colorbar,
        colormap = {seismic}{rgb=(0,0,0.3) rgb=(0,0,1) color=(white) rgb=(1,0,0) rgb=(0.5,0,0)},
        point meta min=-6.099938571800616, 
        point meta max=6.099938571800616,  
        colorbar style={
            ylabel={$\textnormal{sgn}(a)\log(1+|a|)$},
            ytick = {-6,-4,-2,0,2,4,6},
            yticklabel style={
            text width=1em,
            align=right
        }
        }
    ]
    
    \addplot [draw=none] coordinates {(0,0)};

    \end{axis}

\node[] () [below = 1.5em] at (3,-0.5) {}; %

\end{tikzpicture}
    \caption{Experimental results for the mixing time of 25 data sets. For every data set, the shown curve is the average perturbation degree of the Hypercurveball method minus the average perturbation degree of the hyperedge-shuffle method. Positive curves reflect data sets where the Hypercurveball method mixes faster, and negative curves reflect data sets where the hyperedge-shuffle method mixes faster. The x-axis represents the number of steps normalized by the total number of simulated steps for each dataset, allowing them to be displayed on the same scale. The color bar uses $a = \mathds{E}_V[\min\{d_1,d_2\}] - \mathds{E}_E[\min\{d_1,d_2\}]$ and is positive iff $a > 0$.}
    \label{fig:MP_all_data}
\end{figure}
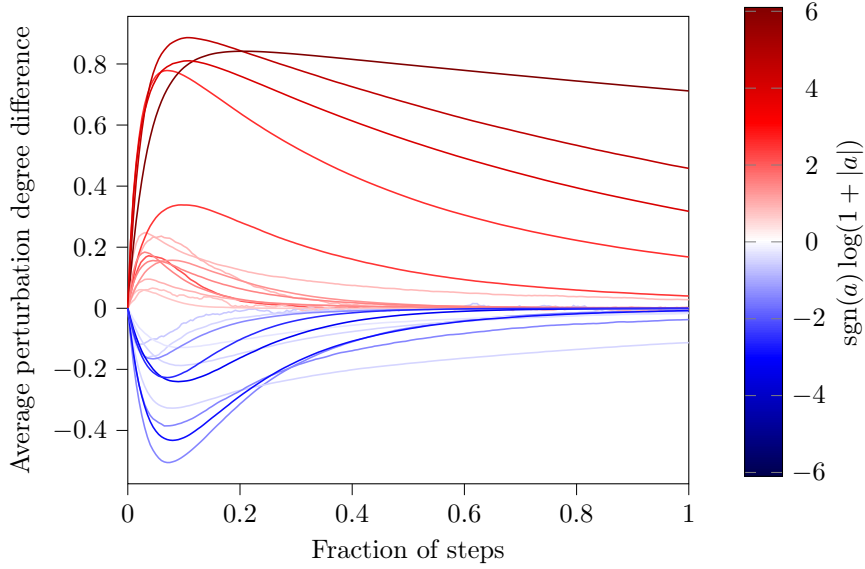

Figure \ref{fig:experiment_artificial} shows that the Hypercurveball method mixes faster on the first artificial data set, and the hyperedge-shuffle method mixes faster on the second and third artificial data set. Therefore, we hypothesize that the Hypercurveball method mixes faster when the expected minimum degree of two nodes is bigger than the expected minimum degree of two hyperedges. For validation of this hypothesis, we perform the numerical experiments on 17 more data sets \cite{data_ceo_club, data_elite, data_kidnappings, data_plant_pol_kato, data_crime, data_plant_pol_robertson, data_NYC_restaurant_checkin_tips, data_contact_primary_school, data_senate_bills, data_house_committees, data_cat_edge_music_blues_reviews, data_madison_vegas_algebra_geometry, data_uchoice_walmart_items}. See Appendix  \ref{app:more_datasets} for the statistics and details of the experiments. Figure \ref{fig:MP_all_data} shows the results for all data sets. From this figure, we conclude that the hypothesis holds for all these data sets. All data used for this analysis can be found here: \cite{algorithms}.

Lastly, we provide an estimate of the scaling of the mixing time for both the Hypercurveball method and the hyperedge-shuffle method by plotting the mixing time of each experiment against a function of the node or hyperedge degrees. Let
\begin{align*}
    f_{\min}^V(\vb*{d}) := \frac{|V|\Bar{d_V}}{2\mathds{E}_V[\min(d_1,d_2)]}\\
    f_{\min}^E(\vb*{d}) := \frac{|E|\Bar{d_E}}{2\mathds{E}_E[\min(d_1,d_2)]},
\end{align*}
where $\Bar{d_V}$ is the average node degree and $\Bar{d_E}$ is the average hyperedge degree. The numerator of $f_{\min}^V(\vb*{d})$, $|V|\Bar{d_V}$, is the sum of the number of elements of all incidence sets, and the denominator $2 \mathds{E}_V[\min(d_1,d_2)]$ is the expected maximum number of elements in incidence sets that can be modified in every step of the algorithm. Then, $f_{\min}^V(\vb*{d})$ can be interpreted as the minimal number of hypertrades needed to modify every entry in every incidence set of the hypergraph. For directed hypergraphs, the results for the tail and head are averaged to obtain one value for $f_{\min}^V(\vb*{d})$. Similarly, $f_{\min}^E(\vb*{d})$ can be interpreted as the minimal number of hyperedge-shuffles needed to modify every entry in every hyperedge set of the hypergraph. Note that $|V|\Bar{d}_V=|E|\Bar{d}_E$.

Since the two methods are symmetric when there are no restrictions on the hyperedge types $d, s$ and $m$, the mixing time of the Hypercurveball method w.r.t. $f_{\min}^V(\vb*{d})$ and that of the hyperedge-shuffle method w.r.t. $f_{\min}^E(\vb*{d})$ are expected the scale similarly. The exact hyperedge connections might influence the mixing time of both methods differently, but here we investigate only the scaling in the global network parameters. Therefore, the Hypercurveball and hyperedge-shuffle results are analyzed in the same plot, where the Hypercurveball data points are plotted against $f_{\min}^V(\vb*{d})$ and the hyperedge-shuffle data points are plotted against $f_{\min}^E(\vb*{d})$. For completeness, we also plotted the extrapolated mixing times for the experiments that did not show mixing yet, by performing a least squares fit of the function $f(x)= L- a \exp(-bx) - c \exp (-dx)$ to the average perturbation degree curve and computing the mixing time of this function. The function seems to fit well, as the Root Mean Square Error of each fit is two to three orders of magnitude smaller than the data. Details on this extrapolation, as well as the measured and extrapolated mixing times are shown in Appendix \ref{app:fit}. See Figure \ref{fig:mixing_scaling} for the results. In this figure, again we see that the neither of the two methods mixes consistently faster than the other. In addition, observe that the mixing time seems to scale linearly on the log-log plot. Therefore, the function $f(x) = ax+b$ is fitted using least squares to the measured data in the log-log plot, which yields $a \approx 1.07 , b \approx 1.70$. This suggests that
\begin{align*}
    \textnormal{Mixing time Hypercurveball} &\approx 5.46 (f_{\min}^V(\vb*{d}))^{1.07}\\
    \textnormal{Mixing time hyperedge-shuffle} &\approx 5.46 (f_{\min}^E(\vb*{d}))^{1.07}.
\end{align*}

\begin{figure}[tbp]
    \centering
    \begin{tikzpicture}

\definecolor{darkgray176}{RGB}{176,176,176}
\definecolor{darkorange25512714}{RGB}{255,127,14}
\definecolor{lightgray204}{RGB}{204,204,204}
\definecolor{steelblue31119180}{RGB}{31,119,180}

\begin{axis}[
width=0.5*\linewidth,
log basis x={10},
log basis y={10},
tick align=outside,
tick pos=left,
x grid style={darkgray176},
xlabel={$f^V_{\min}(\vb*{d})$ or $f^E_{\min}(\vb*{d})$},
xmin=7.98974939662335, xmax=46346.9941389074,
xmode=log,
xtick style={color=black},
y grid style={darkgray176},
ylabel={Mixing time},
ymin=37.4617869425372, ymax=562415.883940937,
ymode=log,
ytick style={color=black}
]
\addplot [draw=steelblue31119180, fill=steelblue31119180, mark=*, only marks]
table{%
x  y
15.3887688984881 87
13.1907433380084 81
182.388888888889 1988
2247.03711455581 23442
268.711133767809 3136
16.1722090261283 103
25.9375 201
981.391369863014 7575
481.766206896552 3736
1490.10083342755 12027
559.932426815051 4519
668.006519554708 8683
1105.97235662722 7881
570.886053589183 4630
852.664507964373 7204
1132.44277569392 9329
449.155413436724 4107
832.933896315768 7939
586.381133911915 6264
161.104501383996 1118
252.712851718293 2865
258.778491738181 1806
};
\addplot [draw=steelblue31119180, fill=steelblue31119180, mark=+, only marks]
table{%
x  y
544.498269991213 6860
3571.16318570491 77877
2187.16169118846 198835
};
\addplot [draw=darkorange25512714, fill=darkorange25512714, mark=*, only marks]
table{%
x  y
11.8467933491686 58
18.9224376731302 110
112.663997536187 1673
139.083305170666 1459
2886.43594360926 24996
2982.56426660077 32220
32.2983870967742 229
45.4184549356223 328
500.99 3915
250.98 1952
130.760607144454 970
418.983540466922 3446
262.82427376298 1782
1375.83451452697 13879
676.851553803787 5247
560.170044324158 4418
1470.04067077636 13141
1097.51477771659 9065
1300.63343859086 9696
};
\addplot [draw=darkorange25512714, fill=darkorange25512714, mark=+, only marks]
table{%
x  y
580.565140974902 21665
17886 314051
2188.93582420478 20404
7136.47631248489 63781
10347.3870138175 78734
31257.4768160894 178073
};
\addplot [semithick, black]
table {%
11.8467933491686 77.4801991335578
649.512712180603 5689.5856917948
1287.17863101204 11852.4792598385
1924.84454984347 18252.3578946919
2562.5104686749 24811.7616762551
3200.17638750634 31492.6051548961
3837.84230633777 38271.9430884179
4475.50822516921 45134.3537492902
5113.17414400064 52068.7250871896
5750.84006283207 59066.6526544908
6388.50598166351 66121.5480516189
7026.17190049494 73228.1020316408
7663.83781932638 80381.9413043253
8301.50373815781 87579.3986421044
8939.16965698924 94817.3529523909
9576.83557582068 102093.11449976
10214.5014946521 109404.340346765
10852.1674134835 116748.970652306
11489.833332315 124125.17975053
12127.4992511464 131531.33794532
12765.1651699778 138965.981229817
13402.8310888093 146427.786971071
14040.4970076407 153915.554155373
14678.1629264721 161428.187169556
15315.8288453036 168964.682358571
15953.494764135 176524.116787789
16591.1606829665 184105.638774442
17228.8266017979 191708.459852129
17866.4925206293 199331.847906306
18504.1584394608 206975.121274284
19141.8243582922 214637.643645526
19779.4902771236 222318.819630593
20417.1561959551 230018.090892258
21054.8221147865 237734.932752094
21692.4880336179 245468.851201444
22330.1539524494 253219.380258074
22967.8198712808 260986.079619785
23605.4857901122 268768.532574283
24243.1517089437 276566.344131142
24880.8176277751 284379.139347023
25518.4835466065 292206.561819713
26156.149465438 300048.272330159
26793.8153842694 307903.9476147
27431.4813031008 315773.279252195
28069.1472219323 323655.972652886
28706.8131407637 331551.746137577
29344.4790595951 339460.330097244
29982.1449784266 347381.466224441
30619.810897258 355314.906808981
31257.4768160894 363260.414091277
};
\end{axis}

\end{tikzpicture}
    \begin{tikzpicture} 

\definecolor{darkgray176}{RGB}{176,176,176}
\definecolor{darkorange25512714}{RGB}{255,127,14}
\definecolor{lightgray204}{RGB}{204,204,204}
\definecolor{steelblue31119180}{RGB}{31,119,180}

    \begin{axis}[%
    hide axis,
    xmin=10,
    xmax=50,
    ymin=0,
    ymax=0.4,
    legend style={draw=white!15!black,legend cell align=left}
    ]
    \addlegendimage{steelblue31119180,mark=*, only marks}
    \addlegendentry{Hypercurveball, measured};
    \addlegendimage{steelblue31119180,mark=+, only marks}
    \addlegendentry{Hypercurveball, extrapolated};
    \addlegendimage{darkorange25512714,mark=*, only marks}
    \addlegendentry{Hyperedge-shuffle, measured};
    \addlegendimage{darkorange25512714,mark=+, only marks}
    \addlegendentry{Hyperedge-shuffle, extrapolated};
    \addlegendimage{semithick, black}
    \addlegendentry{Least squares fit of measured data};
    \end{axis}

\node[] () [below = 1.5em] at (3,0.15) {}; %

\end{tikzpicture}
    \caption{Scaling of the mixing times. The Hypercurveball data points are plotted against $f_{\min}^V(\vb*{d})$ and the hyperedge-shuffle data points are plotted against $f_{\min}^E(\vb*{d})$. Fit: $f(x)=5.46 x^ {1.07}$.}
    \label{fig:mixing_scaling}
\end{figure}
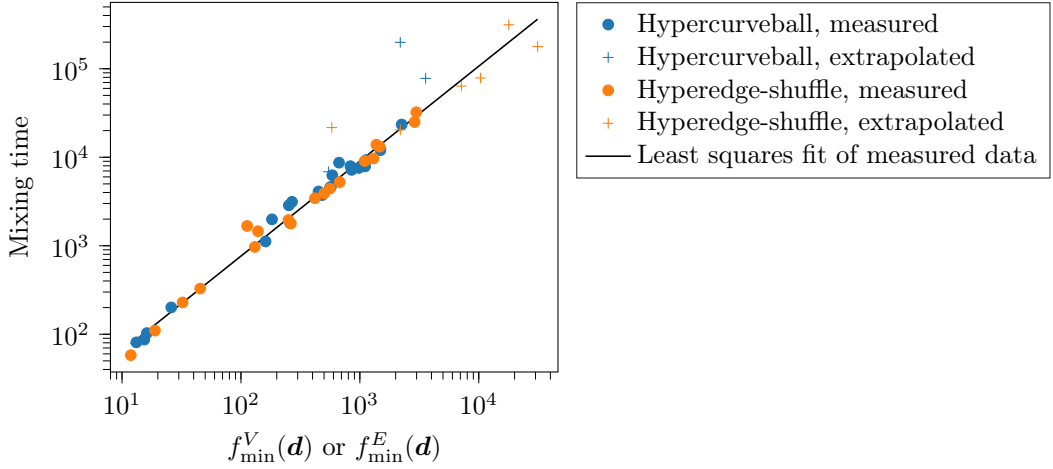

\newpage
\section{Conclusion}
\label{section:conclusion}
In this work, we have introduced the Hypercurveball algorithm as a method to generate randomized (directed) hypergraphs with a fixed degree sequence. This sequential algorithm exchanges elements in two node incidence sets in each step. We show that the Hypercurveball algorithm samples uniformly for some hypergraph spaces, and with bias for others. In particular, the algorithm always samples uniformly when multi-hyperedges are allowed, and when self-loops and multi-hyperedges are allowed. On the other hand, the algorithm may sample with bias from when degenerate hyperedges are allowed and multi-hyperedges are not, and when degenerate hyperedges are allowed and self-loops and/or multi-hyperedges are not. For one undirected and three directed hypergraph spaces the question whether it samples uniformly remains open. Combining existing results with our new results, the graph cases are all solved: for undirected graphs, the algorithm always samples uniformly, except if self-loops are allowed and multi-edges are not. For digraphs, the algorithm always samples uniformly. 

We experimentally compared the mixing time of the Hypercurveball algorithm to that of the hyperedge-shuffle algorithm, which exchanges elements in two hyperedges in each step. Our experiments showed that the Hypercurveball algorithm sometimes mixes faster and sometimes mixes slower than the hyperedge-shuffle algorithm. This contrasts with the results for graphs where the Curveball algorithm typically mixes faster than the edge-swap algorithm. Based on the experiments, we conjecture that the Hypercurveball algorithm mixes faster than the hyperedge-shuffle algorithm if the expected minimum degree of two nodes is bigger than the expected minimum degree of two hyperedges. Lastly, we experimentally show polynomial scaling of the mixing times of both the Hypercurveball and the hyperedge-shuffle algorithm.

\section{Discussion}
\label{section:discussion}

We now provide some discussion of our results.

For four hypergraph spaces, it is still an open question whether the Hypercurveball method always samples uniformly (see Tables \ref{tab:results_undir} and \ref{tab:results_dir}). Since there exist no small counterexamples (and usually the counterexamples are small) for each of these spaces, we conjecture that the method samples uniformly. To prove this, one needs to show that any hypergraph in such a space can be transformed into any other hypergraph in the space using hypertrades, without creating forbidden hyperedge types ($d,s$, and/or, $m$). Unfortunately, related results for graphs heavily use the adjacency matrix to prove such statements. For hypergraphs, the adjacency matrix is not defined (or is a tensor), which makes it necessary to apply new techniques to pursue this.

Furthermore, while we only analyze hypertrades containing two nodes, we believe that more complex type of hypertrade that uses three nodes, could be more powerful. Such a 3-hypertrade would take the incidence sets of three nodes and randomly repartition their union into three sets of the correct sizes, similar to line \ref{line:partition} in Algorithm \ref{alg:trade_d}. For the edge-shuffle algorithm on graphs, a shuffle with three edges was analyzed in \cite{berger2010}. Berger and Müller-Hannemann show that such a 3-shuffle, which they call a \textit{3-cycle reorientation}, allows for uniform sampling of simple (di)graphs. We believe that a 3-hypertrade will yield uniform samples in all hypergraph spaces, and it would be interesting to prove this. On the other hand, such a hypertrade would be computationally more expensive, so optimizing when to use such computationally expensive trades would be another line of further research. Similarly, for the hyperedge-shuffle, a generalization of the 3-cycle reorientation would be more powerful and more computationally expensive than the hyperedge-shuffle using two hyperedges.

Finally, as extensively discussed by \cite{fosdick2018}, dyadic graphs can be either vertex-labeled or stub-labeled. This holds for hypergraphs as well, and we refer to \cite{kraakman2024} for definitions and examples. In this work, we analyzed stub-labeled hypergraphs. As in the work of \cite{kraakman2024}, the results easily extend to vertex-labeled hypergraphs. The hyperedges are considered not labeled.

\section{Proofs}
\label{section:proofs}

\subsection{Stationary distribution of a Markov chain}
The Hypercurveball algorithm can be interpreted as random walk on a Markov chain, where the states are hypergraphs in the applicable hypergraph space and the transitions are given by the hypertrades. The algorithm samples uniformly if the stationary distribution of the Markov chain is uniform, which is the case if the Markov chain is \cite{levin2009}
\begin{enumerate}
    \item regular:\\
    $\sum_{S'} P(S \rightarrow S')$ is equal for all states $S$ and $\sum_{S} P(S \rightarrow S')$ is equal for all states $S'$, where $P(A \rightarrow B)$ is the transition probability between states $A$ and $B$;
    \item aperiodic:\\
    The cycles in the Markov chain have greatest common divisor 1;
    \item strongly connected:\\
    For every two states $S,S'$ there exists a path from $S$ to $S'$ in the Markov chain.
\end{enumerate}

In this section, we prove regularity and aperiodicity of the Markov chains behind the Hypercurveball algorithms for all undirected and directed hypergraph spaces $\mathcal{H}_{x}(\vb*{d})$, $x \subseteq \{s,d,m\}$. To do this, we analyze the graph of graphs $\mathcal{G}(\mathcal{H}_x(\vb*{d}))$, which denotes the transition graph underlying the Markov chain on the space $\mathcal{H}_x(\vb*{d})$. In this graph, the vertices represent the states and the arcs represent the transitions of the Markov chain.
\\
\begin{lemma}
\label{lemma:regular}
    The Markov chains behind the Hypercurveball algorithms for the undirected and directed hypergraph spaces $\mathcal{H}_{x}(\vb*{d})$, $x \subseteq \{s,d,m\}$, are regular.
\end{lemma}
\begin{proof}
    Let $x \subseteq \{s,d,m\}$. By construction of the algorithm,
    \begin{align}
    \label{eq:sum_out_=1}
        \sum_{S' \in \mathcal{H}_{x} (\vb*{d})} P(S \rightarrow S') = 1
    \end{align}
    for all $S \in \mathcal{H}_{x} (\vb*{d})$. In addition, if $P(S \rightarrow S') > 0$ and $S \neq S'$ then to transition from $S$ to $S'$ using one hypertrade, the two nodes whose incidence sets differ between $S$ and $S'$ should be selected in the hypertrade. Then, during the hypertrade, the combined incidence set should be partitioned into the correct sets for both nodes. This results in    
    \begin{align*}
        P(S \rightarrow S') &= \binom{|V(S)|}{2}^{-1} \binom{d_v + d_w}{d_v} ^{-1}\\
        &= \binom{|V(S')|}{2}^{-1} \binom{d_{v'} + d_{w'}}{d_{v'}} ^{-1}\\
        &= P(S'\rightarrow S).
    \end{align*}
    For \texttt{hypertrade\_nodeg} and \texttt{hypertrade\_simple}, these transition probabilities are slightly different, since $\Tilde{I} \cup \Tilde{J}$ is used for partitioning. However, the size of this set for the transition $S \rightarrow S'$ is equal to the size of this set for the reversed transition $S' \rightarrow S$, so the equalities still hold. Therefore, 
    \begin{align*}
        \sum_{S \in \mathcal{H}_{x} (\vb*{d})} P(S \rightarrow S') = \sum_{S \in \mathcal{H}_{x} (\vb*{d})} P(S'\rightarrow S) = 1,
    \end{align*}
    by (\ref{eq:sum_out_=1}). For a directed hypergraph space, the proof can be adjusted by using 
    \begin{align*}
        P(S \rightarrow S') &= \binom{|V(S)|}{2}^{-1} \binom{d_v^{\tail} + d_w^{\tail}}{d_v^{\tail}} ^{-1} \binom{d_v^{\head} + d_w^{\head}}{d_v^{\head}} ^{-1}.
    \end{align*}
\end{proof}

\begin{lemma}
\label{lemma:aperiodic}
    The Markov chains behind the Hypercurveball algorithms for the undirected and directed hypergraph spaces $\mathcal{H}_{x}(\vb*{d})$, $x \subseteq \{s,d,m\}$, are aperiodic.
\end{lemma}
\begin{proof}
    Let $x \subseteq \{s,d,m\}$. By construction of the algorithm, $P(S \rightarrow S)>0$ for all $S \in \mathcal{H}_x(\vb*{d})$. Therefore, the graph of graphs $\mathcal{G}(\mathcal{H}_x(\vb*{d}))$ contains self-loops, which makes it aperiodic.
\end{proof}

In the remainder of Section \ref{section:proofs}, we prove or disprove that certain Markov chains are strongly connected. This finalizes the proofs of Theorems \ref{thm:uniform_sampling}, \ref{thm:uniform_sampling_graph}, \ref{thm:biased_sampling} and \ref{thm:methods_not_equivalent}.

\subsection{Proof of uniformity: Theorem \ref{thm:uniform_sampling}}

\subsection{Proof of uniformity for graphs: Theorem \ref{thm:uniform_sampling_graph}}
For the proof of Theorem \ref{thm:uniform_sampling_graph}, we first note that an edge swap in either an undirected or a directed graph is also a simple hyperedge shuffle, as in Definition \ref{def:simple}. By Lemma \ref{lemma:simple}, the result of an edge-swap can be achieved using a Curveball trade. 

\begin{figure}[tbp]
\centering
\begin{subfigure}{0.43\textwidth}
    \centering
    \newcommand{\offsetx}{4*\tikzscale}
    
\begin{tikzpicture}[-Stealth, line width=1.3pt,auto,
                    thick,main node/.style={circle,draw, minimum size=\tikznodesize}, inner sep=1pt]

  \node[main node] (a) at (0,0) {$a$};
  \node[main node] (b) at (2*\tikzscale,0) {$b$};
  \node[main node] (c) at (\tikzscale,1.73*\tikzscale) {$c$};

  \path[every node/.style={font=\sffamily\small}]
    (a) edge node {} (b)
    (b) edge node {} (c)
    (c) edge node {} (a);

\end{tikzpicture}
    \caption{A directed triangle, orientation anti-clockwise.}
    \label{fig:triangle_anticlockwise}
\end{subfigure} \hspace{0.1\textwidth}
\begin{subfigure}{0.43\textwidth}
    \centering
    \newcommand{\offsetx}{4*\tikzscale}
    
\begin{tikzpicture}[-Stealth, line width=1.3pt,auto,
                    thick,main node/.style={circle,draw, minimum size=\tikznodesize}, inner sep=1pt]

  \node[main node] (a) at (0,0) {$a$};
  \node[main node] (b) at (2*\tikzscale,0) {$b$};
  \node[main node] (c) at (\tikzscale,1.73*\tikzscale) {$c$};

  \path[every node/.style={font=\sffamily\small}]
    (a) edge node {} (c)
    (b) edge node {} (a)
    (c) edge node {} (b);

\end{tikzpicture}
    \caption{A directed triangle, orientation clockwise.}
    \label{fig:triangle_clockwise}
\end{subfigure}
\caption{The two realizations of the digraph degree sequence $\vb*{d}=(\vb*{d}_V,\vb*{d}_E)$ with $\vb*{d}_V=((2,2),(2,2),(2,2))$ in the digraph space without self-loops.}
\label{fig:example_triangles}
\end{figure}
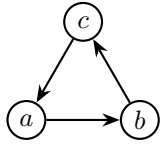
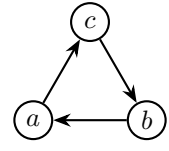

The Curveball algorithm for graphs allows for more graph transformations than the edge-swap algorithm. In particular, every edge-swap can be achieved using a Curveball trade (see Lemma \ref{lemma:simple}). Moreover, the triangle swap in digraphs, which reverses the direction of a directed triangle, cannot be achieved using a single edge-swap, but can be achieved using a single Curveball trade. 
\\
\begin{lemma}
\label{lemma:triangle_swap}
    Let a triangle swap be a transition of three edges forming a directed triangle in one orientation ($a \rightarrow b$, $b \rightarrow c$ and $c \rightarrow a$) to the directed triangle in reversed orientation ($a \rightarrow c$, $b \rightarrow a$ and $c \rightarrow b$). An edge shuffle cannot achieve a triangle swap, and a Curveball trade can. 
\end{lemma}
\begin{proof}
    Let $H$ be the directed triangle as in Figure \ref{fig:triangle_anticlockwise}. Performing an edge shuffle on the directed edges $c \rightarrow a$ and $a \rightarrow b$ results in an unchanged digraph or in a digraph with a self-loop $a \rightarrow a$. Similarly, the other two edge pairs result in an unchanged digraph or a digraph with a self-loop. Therefore, one edge shuffle cannot transform the triangle in Figure \ref{fig:triangle_anticlockwise} into the triangle in Figure \ref{fig:triangle_clockwise}. A Curveball trade, however, can do this. For example, let $e_1 = (c,a), e_2=(a,b), e_3=(b,c)$. Then, $I_a(H)=(e_2,e_1), I_b(H)=(e_3,e_2), I_c(H)=(e_1,e_3)$. Performing a trade on the nodes $a$ and $b$ can result in the hypergraph $H'$ with incidence sets $I_a(H') = (e_3,e_2), I_b(H')=(e_2,e_1)$, i.e., the directed triangle in Figure \ref{fig:triangle_clockwise}.
\end{proof}

Now, we present a result on the digraph space $\mathcal{H}_m(\vb*{d}'_2)$. After, we present the proof of Theorem \ref{thm:uniform_sampling_graph}.
\\
\begin{lemma}
\label{lemma:strongly_connected_H_m}
    The Markov chain underlying the Curveball algorithm on the digraph space $\mathcal{H}_m(\vb*{d}'_2)$ is strongly connected.
\end{lemma}
\begin{proof}
    Let $G_1,G_2 \in \mathcal{H}_{m}(\vb*{d}'_2)$. We will prove that $G_1$ can be transformed into $G_2$ using Curveball trades, without creating digraphs with self-loops. To that end, we analyze the symmetric difference $G_1 \Delta G_2 := (E(G_1)\backslash E(G_2)) \cup (E(G_2) \backslash E(G_1))$ and perform a proof by induction on the number of edges in $G_1 \Delta G_2$. We use that the symmetric difference consists of alternating cycles: $(v_1,v_2,\hdots,v_{2n-1},v_1)$, where $(v_i,v_{i+1}) \in E(G_1), (v_{i-1},v_{i}) \in E(G_2)$ for $i$ even \cite{berger2010}. 
    
    Base case: $|G_1 \Delta G_2| = 4$. Then, $G_1 \Delta G_2 = \{(a,b),(c,b),(c,d),(a,d)\}$, where $(a,b),(c,d)\in E(G_1)\backslash E(G_2)$ and $(c,b),(a,d) \in E(G_2) \backslash E(G_1)$. Let $e_1=(a,b), e_2=(c,d)$. A Curveball trade on $G_1$ with nodes $a$ and $c$ can result in graph $G_2$, since
    \begin{align*}
        I_a(G_2) &= (I^{\tail}_a(G_1) \backslash \{e_1\} \uplus \{e_2\}, I^{\head}_a(G_1))\\
        I_c(G_2) &= (I^{\tail}_c(G_1) \backslash \{e_2\} \uplus \{e_1\}, I^{\head}_c(G_1)).
    \end{align*}

    Induction hypothesis: Any $G_1$ can be transformed into $G_2$ using a sequence of Curveball trades if $|G_1 \Delta G_2| \leq 2k$, with $k \geq 2$. We prove that $G_1$ can be transformed into $G_2$ using a sequence of Curveball trades if $|G_1 \Delta G_2| = 2(k+1)$. By construction, if $G_1 \neq G_2$ then $G_1 \Delta G_2$ contains at least 3 nodes. If $G_1 \Delta G_2$ contains exactly 3 nodes, then it looks as illustrated in Figure \ref{fig:sym_dif_3_nodes}, i.e., $G_1$ contains a triangle with multiplicity $m$, and $G_2$ contains this triangle with reversed direction with multiplicity $m$. Then, one Curveball trade using nodes $b$ and $c$ can transform $G_1$ into $G_2$, as we showed for a single triangle in Lemma \ref{lemma:triangle_swap}. Now, assume $G_1 \Delta G_2$ consists of at least 4 nodes. If $G_1 \Delta G_2$ consists of multiple components $C_1,C_2,\hdots,C_n$, with $n \geq 2$, then let $\Tilde{G}$ be the digraph such that $G_1 \Delta \Tilde{G} = G_1 \Delta G_2 \backslash C_1$ and $\Tilde{G} \Delta G_2 = C_1$. This graph $\tilde{G}$ belongs to $\mathcal{H}_m(\vb*{d})$ as well, since $C_1$ consists of alternating cycles \cite{kraakman2024}. Now, $|G_1 \Delta \Tilde{G}| \leq 2k$ and $|\Tilde{G} \Delta G_2| \leq 2k$, so by the induction hypothesis $G_1$ can be transformed into $\Tilde{G}$ and $\Tilde{G}$ can be transformed into $G_2$.     

    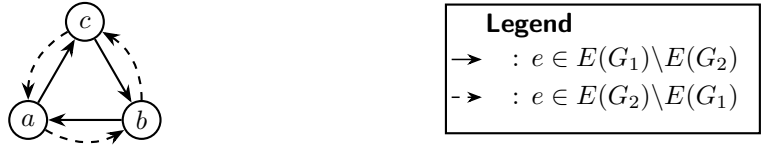
\begin{figure}[tb]
        \centering
        \newcommand{\offsetx}{1}
\newcommand{\offsety}{0.5}
    
\begin{tikzpicture}[-Stealth, line width=1.3pt,auto,
                    thick,main node/.style={circle,draw, minimum size=\tikznodesize}, inner sep=1pt]

  \node[main node] (a) at (0,0) {$a$};
  \node[main node] (b) at (2*\tikzscale,0) {$b$};
  \node[main node] (c) at (\tikzscale,1.73*\tikzscale) {$c$};

  \path[every node/.style={font=\sffamily\small}]
    (a) edge[dashed, bend right] node {} (b)
    (b) edge[dashed, bend right] node {} (c)
    (c) edge[dashed, bend right] node {} (a)
    (a) edge[] node {} (c)
    (c) edge[] node {} (b)
    (b) edge[] node {} (a);

  \node[anchor=south west, font=\bfseries\sffamily] at (5.0 + \offsetx,0.55+\offsety) {Legend};

  \node at (6.9 + \offsetx,0.3 + \offsety) {: \textcolor{black}{$e \in E(G_1) \backslash E(G_2)$}};
  \draw[line width=0.02cm] (4.6 + \offsetx,0.3 + \offsety) -- (5 + \offsetx,0.3 + \offsety);
  \node at (6.9 + \offsetx,-0.2 + \offsety) {: \textcolor{black}{$e \in E(G_2) \backslash E(G_1)$}};
  \draw[dashed, line width=0.02cm] (4.6 + \offsetx,-0.2 + \offsety) -- (5 + \offsetx,-0.2 + \offsety);

  \node[draw, rectangle, fit={(4.7 + \offsetx,0.85 + \offsety) (8.5 + \offsetx,-0.5 + \offsety)}, inner sep=5pt] {};

\end{tikzpicture}
        \caption{The symmetric difference $G_1 \Delta G_2$ when it contains exactly 3 nodes. Every edge has multiplicity $m$, for some $m \in \mathds{N}$.}
        \label{fig:sym_dif_3_nodes}
    \end{figure}
    
    If $G_1 \Delta G_2$ consists of exactly one component which contains at least four nodes, then there exists a node-disjoint path $P=(a,b,c,d)$ in $G_1 \Delta G_2$ with $(a,b),(c,d)\in E(G_1)\backslash E(G_2)$ and $(c,b) \in E(G_2) \backslash E(G_1)$, or vice versa. Indeed, if this would not be the case, then every alternating path of size 3 is of the form $(a,b,c,a)$, not allowing a fourth node in the connected component. This is a contradiction, thus, the path $P$ must exist. W.l.o.g., let $(a,b),(c,d)\in E(G_1)\backslash E(G_2)$ and $(c,b) \in E(G_2) \backslash E(G_1)$. Then, let $\Tilde{G}$ be the graph with $E(\tilde{G})=E(G_1)\backslash \{(a,b),(c,d)\} \uplus \{(c,b),(a,d)\}$. Then, 
    \begin{align*}
        |G_1 \Delta \tilde{G}| &\leq |(G_1 \Delta G_2) \backslash \{(a,b),(c,d),(c,b)\} \uplus \{(a,d)\}| = 2(k+1)-2 = 2k,\\
        |\tilde{G} \Delta G_2| &= |\{(a,b),(c,d),(c,b),(a,d)\}| = 4 \leq 2k.
    \end{align*}
    Therefore, $G_1$ can be transformed into $\tilde{G}$ and $\tilde{G}$ can be transformed into $G_2$, concluding the induction proof.
\end{proof}

\begin{proof}[Proof of Theorem \ref{thm:uniform_sampling_graph}]
    The undirected cases $\mathcal{H}_{d,m}(\vb*{d}'_1)$ and $\mathcal{H}_{m}(\vb*{d}'_1)$ and the directed case $\mathcal{H}_{s,m}(\vb*{d}'_2)$ follow directly from Theorem \ref{thm:uniform_sampling}. The proof for the undirected space $\mathcal{H}(\vb*{d}'_1)$ is given by \cite{carstens2016}, and the proof for the directed space $\mathcal{H}_s(\vb*{d}'_2)$ is given by \cite{carstens2015}. Here, we present the proofs for the directed spaces $\mathcal{H}_m(\vb*{d}'_2)$ and $\mathcal{H}(\vb*{d}'_2)$.
    
    By Lemmas \ref{lemma:regular} and \ref{lemma:aperiodic}, the Markov chains underlying the Curveball algorithm for the directed spaces $\mathcal{H}_{m}(\vb*{d}'_2)$ and $\mathcal{H}_{}(\vb*{d}'_2)$ are regular and aperiodic. By Lemma \ref{lemma:strongly_connected_H_m}, the Markov chain underlying the Curveball algorithm for the digraph $\mathcal{H}_{m}(\vb*{d}'_2)$ is strongly connected. Therefore, the Curveball method samples uniformly from this space. It remains to prove that the digraph space $\mathcal{H}_{}(\vb*{d}'_2)$ is strongly connected via Curveball trades.

    For the digraph space $\mathcal{H}_{}(\vb*{d}'_2)$, the Markov chain is strongly connected via edge-swaps if triangle-swaps are also allowed \cite{berger2010}. A triangle-swap reverses the direction of a triangle $(a,b),(b,c),(c,a)$ to obtain $(a,c),(c,b),(b,a)$. Such a triangle swap can be achieved using a Curveball trade (Lemma \ref{lemma:triangle_swap}). Combining this with Lemma \ref{lemma:simple}, this proves that the Markov chain is strongly connected through Curveball trades.
\end{proof}

\subsection{Proof of non-uniformity: Theorem \ref{thm:biased_sampling}}
We show that there exist degree sequences $\tilde{\vb*{d}},\hat{\vb*{d}},\vb*{d}^*$ such that the Hypercurveball algorithms sample with bias from the undirected hypergraph space $\mathcal{H}_{d}(\tilde{\vb*{d}})$ and the directed hypergraph spaces $\mathcal{H}_{s,d}(\hat{\vb*{d}}), \mathcal{H}_{d,m}(\vb*{d}^*),\mathcal{H}_{d}(\vb*{d}^*)$.

\begin{proof}[Proof of Theorem \ref{thm:biased_sampling}]
    For the undirected hypergraph space $\mathcal{H}_{d}(\tilde{\vb*{d}})$, consider the degree sequence $\tilde{\vb*{d}} = (\vb*{d}_V,\vb*{d}_E)$, where $\vb*{d}_V = (2,2,2)$ and $\vb*{d}_E = (2,2,2)$. Let $H_1 = (V,E_1)$, $H_2 = (V,E_2)$, with $V=\{u,v,w\}$ and
    \begin{align*}
        E_1&=\{\{u,u\},\{v,v\}, \{w,w\}\}\\
        E_2&=\{\{u,v\},\{v,w\},\{u,w\}\}.
    \end{align*}
    The two hypergraphs $H_1$ and $H_2$ are shown in Figure \ref{fig:biased_undir_H_d}.

    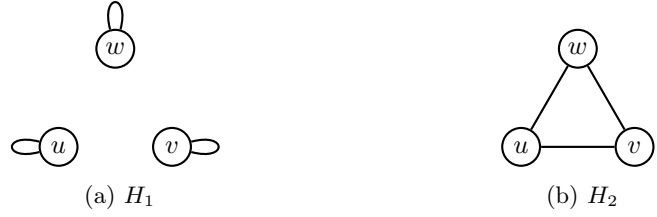
\begin{figure}[tbp]
    \centering
    \begin{subfigure}{0.4\textwidth}
    \centering
        \newcommand{\offsetx}{1}
\newcommand{\offsety}{0.5}
    
\begin{tikzpicture}[line width=1.3pt,auto,
                    thick,main node/.style={circle,draw, minimum size=\tikznodesize}, inner sep=1pt,every loop/.style={}]

  \node[main node] (a) at (0,0) {$u$};
  \node[main node] (b) at (2*\tikzscale,0) {$v$};
  \node[main node] (c) at (\tikzscale,1.73*\tikzscale) {$w$};

  \path[every node/.style={font=\sffamily\small}]
    (a) edge[loop left] node {} (a)
    (b) edge[loop right] node {} (b)
    (c) edge[loop above] node {} (c);
\end{tikzpicture}
        \caption{$H_1$}
    \end{subfigure}
    \begin{subfigure}{0.4\textwidth}
    \centering
        \newcommand{\offsetx}{1}
\newcommand{\offsety}{0.5}
    
\begin{tikzpicture}[line width=1.3pt,auto,
                    thick,main node/.style={circle,draw, minimum size=\tikznodesize}, inner sep=1pt,every loop/.style={}]

  \node[main node] (a) at (0,0) {$u$};
  \node[main node] (b) at (2*\tikzscale,0) {$v$};
  \node[main node] (c) at (\tikzscale,1.73*\tikzscale) {$w$};

  \path[every node/.style={font=\sffamily\small}]
    (a) edge node {} (b)
    (b) edge node {} (c)
    (c) edge node {} (a);
\end{tikzpicture}
        \caption{$H_2$}
    \end{subfigure}
    \caption{Two undirected hypergraphs $H_1,H_2 \in \mathcal{H}_d(\tilde{\vb*{d}})$ which are not connected using hypertrades.}
    \label{fig:biased_undir_H_d}
    \end{figure}
    
    We now show that there exists no path from $H_1$ to $H_2$ in $\mathcal{G}(\mathcal{H}_{d}(\tilde{\vb*{d}}))$. To that end, consider any two nodes $m,n \in V$. Then, $I_m(H_1) = \{e_1,e_1\}$ and $I_n(H_1) = \{e_2,e_2\}$, for two hyperedges $e_1,e_2 \in E_1$. A hypertrade on $m,n$ can result in a hypergraph $H'$ with $I_m(H') = \{e_1,e_2\}$ and $I_n(H') = \{e_1,e_2\}$, which is not allowed since $e_1$ and $e_2$ are then multi-edges, or in a hypergraph $H''$ with $I_m(H'') = \{e_2,e_2\}$ and $I_n(H'') = \{e_1,e_1\}$, which is equivalent to $H$. Therefore, $H_1$ is an isolated vertex in $\mathcal{G}(\mathcal{H}_{d}(\tilde{\vb*{d}}))$. In particular, there exists no path from $H_1$ to $H_2$ in $\mathcal{G}(\mathcal{H}_{d}(\tilde{\vb*{d}}))$.

    For the directed hypergraph space $\mathcal{H}_{s,d}(\hat{\vb*{d}})$, consider the degree sequence $\hat{\vb*{d}} = (\vb*{d}_V,\vb*{d}_E)$, where $\vb*{d}_V = ((2,0),(2,0),(2,0),(0,3))$ and $\vb*{d}_E = ((2,1),(2,1),(2,1))$. Let $H_1 = (V,E_1)$ and $H_2=(V,E_2)$, with $V=\{u,v,w,x\}$ and
    \begin{align*}
        E_1&=\{(\{u,u\},\{x\}),(\{v,v\},\{x\}), (\{w,w\},\{x\})\}\\
        E_2&=\{(\{u,v\},\{x\}),(\{v,w\},\{x\}),(\{u,w\},\{x\})\}.
    \end{align*}
    The two hypergraphs $H_1$ and $H_2$ are shown in Figure \ref{fig:biased_dir_H_sd} \cite{kraakman2024}.

    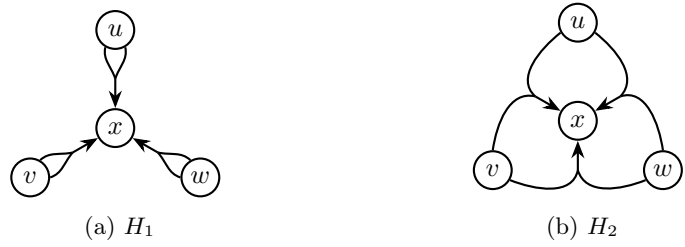
\begin{figure}[tbp]
    \centering
    \begin{subfigure}{0.4\textwidth}
    \centering
      
\begin{tikzpicture}[-Stealth, line width=1.3pt,auto,
                    thick,main node/.style={circle,draw, minimum size=\tikznodesize}, inner sep=1pt]

  \node[main node] (u) at (1.5* \tikzscale, 2.60*\tikzscale) {$u$};
  \node[main node] (v) at (0,0) {$v$};
  \node[main node] (w) at (3*\tikzscale,0) {$w$};
  \node[main node] (x) at (1.5*\tikzscale, 0.87*\tikzscale) {$x$};

\coordinate (uux) at (1.5*\tikzscale,1.74*\tikzscale);

\draw[-,bend left, in=180] (u.300) to (uux);
\draw[-,bend right, in=180] (u.240) to (uux);
\draw[] (uux) to (x);

\coordinate (vvx) at (0.75*\tikzscale,0.44*\tikzscale);

\draw[-,bend left, in=180] (v.60) to (vvx);
\draw[-,bend right, in=180] (v.0) to (vvx);
\draw[] (vvx) to (x);

\coordinate (wwx) at (2.25*\tikzscale,0.44*\tikzscale);

\draw[-,bend left, in=180] (w.180) to (wwx);
\draw[-,bend right, in=180] (w.120) to (wwx);
\draw[] (wwx) to (x);

\end{tikzpicture}
        \caption{$H_1$}
    \end{subfigure}
    \begin{subfigure}{0.4\textwidth}
    \centering
      
\begin{tikzpicture}[-Stealth, line width=1.3pt,auto,
                    thick,main node/.style={circle,draw, minimum size=\tikznodesize}, inner sep=1pt]

  \node[main node] (u) at (1.5* \tikzscale, 2.60*\tikzscale) {$u$};
  \node[main node] (v) at (0,0) {$v$};
  \node[main node] (w) at (3*\tikzscale,0) {$w$};
  \node[main node] (x) at (1.5*\tikzscale, 0.87*\tikzscale) {$x$};

\coordinate (uvx) at (0.75*\tikzscale,1.3*\tikzscale);

\draw[-,bend right, in=250] (u) to (uvx);
\draw[-,bend left, in=115] (v) to (uvx);
\draw[] (uvx) to (x);

\coordinate (uwx) at (2.25*\tikzscale,1.3*\tikzscale);

\draw[-,bend left, in=115] (u) to (uwx);
\draw[-,bend right, in=250] (w) to (uwx);
\draw[] (uwx) to (x);

\coordinate (vwx) at (1.5*\tikzscale,0*\tikzscale);

\draw[-,bend left, in=115] (w) to (vwx);
\draw[-,bend right, in=250] (v) to (vwx);
\draw[] (vwx) to (x);

\end{tikzpicture}
        \caption{$H_2$}
    \end{subfigure}
    \caption{Two directed hypergraphs $H_1,H_2 \in \mathcal{H}_{s,d}(\hat{\vb*{d}})$ which are not connected using hypertrades.}
    \label{fig:biased_dir_H_sd}
    \end{figure}
    
    We now show that there exists no path from $H_1$ to $H_2$ in $\mathcal{G}(\mathcal{H}_{s,d}(\hat{\vb*{d}}))$. To that end, consider any two nodes $m,n \in \{u,v,w\}$. Then $I_m(H_1) = (\{e_1,e_1\},\emptyset), I_n(H_1) = (\{e_2,e_2\},\emptyset)$, for two hyperedges $e_1,e_2 \in E_1$. A hypertrade on $m$ and $n$ can result in a hypergraph $H'$ with $I_m(H') = (\{e_1,e_2\},\emptyset)$ and $I_n(H') = (\{e_1,e_2\},\emptyset)$, which is not allowed since $e_1$ and $e_2$ are then multi-hyperedges, or in a hypergraph $H''$ with $I_m(H'') = (\{e_2,e_2\},\emptyset)$ and $I_n(H'') = (\{e_1,e_1\},\emptyset)$, which is equivalent to $H$. A hypertrade on $x$ and $m$, where $m\in \{u,v,w\}$, results in $H$ as well. Therefore, $H_1$ is an isolated vertex in $\mathcal{G}(\mathcal{H}_{s,d}(\hat{\vb*{d}}))$. In particular, there exists no path from $H_1$ to $H_2$ in $\mathcal{G}(\mathcal{H}_{s,d}(\hat{\vb*{d}}))$.

    For the directed hypergraph spaces $\mathcal{H}_{d,m}(\vb*{d}^*)$ and $\mathcal{H}_{d}(\vb*{d}^*)$, consider the degree sequence $\vb*{d}^* = (\vb*{d}_V,\vb*{d}_E)$, where $\vb*{d}_V = ((2,2),(1,0),(0,1))$ and $\vb*{d}_E = ((2,2),(1,1))$. Let $H_1 = (V,E_1)$ and $H_2=(V,E_2)$, with $V=\{u,v,w\}$ and
    \begin{align*}
        E_1&=\{(\{u,v\},\{u,u\}),(\{u\},\{w\})\}\\
        E_2&=\{(\{u,u\},\{u,w\}),(\{v\},\{u\}).
    \end{align*}
   The two hypergraphs $H_1$ and $H_2$ are shown in Figure \ref{fig:biased_dir_H_dm}.
    
    We now show that there exists no path from $H_1$ to $H_2$ in $\mathcal{G}(\mathcal{H}_{d,m}(\vb*{d}^*))$ or in $\mathcal{G}(\mathcal{H}_{d}(\vb*{d}^*))$. To that end, let $e_1 = (\{u,v\},\{u,u\})$ and $e_2=(\{u\},\{w\})$. We obtain $I_u=(\{e_1,e_2\},\{e_1,e_1\})$, $I_v=(\{e_1\},\emptyset)$ and $I_w=(\emptyset, \{e_2\})$. Now, consider trading $u$ and $v$. This can result in a hypergraph $H'$ with $I_u(H') = (\{e_1,e_1\},\{e_1,e_1\})$, $I_v(H')=(\{e_2\},\emptyset)$, which creates self-loop $e_1$, or in the unchanged hypergraph $H_1$. Now, consider trading $u$ and $w$. This can result in a hypergraph $H''$ with $I_u(H'') = (\{e_1,e_2\},\{e_1,e_2\})$, $I_w(H'')=(\emptyset,\{e_1\})$, which creates self-loop $e_2$, or in the unchanged hypergraph $H_1$. Trading $v$ and $w$ does not change anything. Therefore, $H_1$ is an isolated vertex in $\mathcal{G}(\mathcal{H}_{d,m}(\vb*{d}^*))$ and in $\mathcal{G}(\mathcal{H}_{d}(\vb*{d}^*))$. In particular, there exists no path from $H_1$ to $H_2$ in $\mathcal{G}(\mathcal{H}_{d,m}(\vb*{d}^*))$ or in $\mathcal{G}(\mathcal{H}_{d}(\vb*{d}^*))$.
\end{proof}

\subsection{Proof of non-equivalence: Theorem \ref{thm:methods_not_equivalent}}
\label{section:pf_thm4}
To prove Theorem \ref{thm:methods_not_equivalent}, we present two degree sequences for the directed hypergraph spaces $\mathcal{H}_{d}(\vb*{d}_2)$ and $\mathcal{H}_{d,m}(\vb*{d}_2)$ where one is sampled uniformly with the Hypercurveball method and with bias with the hyperedge-shuffle method, and the other is sampled uniformly with the hyperedge-shuffle method and with bias with the Hypercurveball method. 

\begin{proof}[Proof of Theorem \ref{thm:methods_not_equivalent}]
    First, consider the degree sequence $\vb*{d}^* = (\vb*{d}_V,\vb*{d}_E)$, with \\$\vb*{d}_V=((2,2),(1,0),(0,1))$ and $\vb*{d}_E = ((2,2),(1,1))$. The only two hypergraphs in $\mathcal{H}_{d}(\vb*{d}^*)=\mathcal{H}_{d,m}(\vb*{d}^*)$ are $H_1$ and $H_2$ as described in the proof of Theorem \ref{thm:biased_sampling} and shown in Figure \ref{fig:biased_dir_H_dm}. As shown in the proof of Theorem \ref{thm:biased_sampling}, the Hypercurveball method samples with bias from this space. On the other hand, the two hypergraphs are connected in the Markov chain of the hyperedge-shuffle, as there only are two hyperedges. Therefore, the hyperedge-shuffle method samples uniformly from this space. 

    Next, consider the degree sequence $\vb*{d}^{**} = (\vb*{d}_V,\vb*{d}_E)$, where $\vb*{d}_V = ((3,3),(2,2))$ and $\vb*{d}_E = ((2,2),(2,2),(1,1))$. The only two hypergraphs in $\mathcal{H}_{d}(\vb*{d}^{**})=\mathcal{H}_{d,m}(\vb*{d}^{**})$ are $H_1=(V,E_1)$ and $H_2=(V,E_2)$, where $V=\{u,v\}$ and
    \begin{align*}
        E_1 &= \{(\{v,v\},\{u,u\}),(\{u,u\},\{u,v\}),(\{u\},\{v\})\}\\
        E_2 &= \{(\{u,v\},\{u,u\}),(\{u,u\},\{v,v\}),(\{v\},\{u\})\}.
    \end{align*}
    The hypergraphs $H_1$ and $H_2$ are shown in Figure \ref{fig:noneq_fig2}. The Hypercurveball method samples uniformly from this space, since there are only two nodes. On the other hand, the hyperedge-shuffle method samples with bias from this space. Let $e_1 = (\{v,v\},\{u,u\}), e_2 = (\{u,u\},\{u,v\}), e_3=(\{u\},\{v\})$. A hyperedge-shuffle with $e_1$ and $e_2$ can result in 
    \begin{align*}
        e'_1 &= (\{u,v\},\{u,u\})\\
        e'_2 &= (\{u,v\},\{u,v\}),
    \end{align*}
    which creates self-loop $e'_2$. A hyperedge-shuffle with $e_1$ and $e_3$ can result in 
    \begin{align*}
        e'_1 &= (\{u,v\},\{u,u\})\\
        e'_3 &= (\{v\},\{v\}),
    \end{align*}
    which creates self-loop $e'_3$, or
    \begin{align*}
        e'_1 &= (\{v,v\},\{u,v\})\\
        e'_3 &= (\{u\},\{u\}),
    \end{align*}
    which creates self-loop $e'_3$, or 
    \begin{align*}
        e'_1 &= (\{u,v\},\{u,v\})\\
        e'_3 &= (\{v\},\{u\}),
    \end{align*}
    which creates self-loop $e'_1$. A hyperedge-shuffle with $e_2$ and $e_3$ can result in
    \begin{align*}
        e'_2 &= (\{u,u\},\{v,v\})\\
        e'_3 &= (\{u\},\{u\}),
    \end{align*}
    which creates self-loop $e'_3$. Therefore, $H_1$ is an isolated vertex in the hyperedge-shuffle graphs $\mathcal{G}(\mathcal{H}_{d}(\vb*{d}^{**}))$ and $\mathcal{G}(\mathcal{H}_{d,m}(\vb*{d}^{**}))$. In particular, there is exists no path from $H_1$ to $H_2$ in the hyperedge-shuffle graphs $\mathcal{G}(\mathcal{H}_{d}(\vb*{d}^{**}))$ or $\mathcal{G}(\mathcal{H}_{d,m}(\vb*{d}^{**}))$.
\end{proof}

\paragraph{Funding.} This work was supported by Nederlandse Organisatie voor Wetenschappelijk Onderzoek [M2 grant 0.379].

\newpage

\setcounter{figure}{0}
\renewcommand{\thefigure}{A\arabic{figure}}

\setcounter{table}{0}
\renewcommand{\thetable}{A\arabic{table}}

\appendix

\section{Biased hypertrades}
\label{app:bias}
As explained in Section \ref{section:curveball}, the \texttt{hypertrade} function cannot easily be adjust to not create any self-loops or multi-hyperedges. Here, we explain this in more detail and present examples.

In all the examples, we show that the difficulty lies in creating a uniformly random partition which satisfies constraints. Creating the partition sequentially creates bias. Creating the partition without considering the constraints and checking afterwards whether the constraints are met does not create bias and is similar to the strategy presented in Algorithm \ref{alg:curveball_overig} (line \ref{line:check}).

\subsection{\texttt{hypertrade} not creating multi-edges in undirected hypergraphs in \texorpdfstring{$\mathcal{H}_{s,d}(\vb*{d}_1)$}{H-sd}}
Consider the example where $I_v \uplus I_w = \{x,x,y,a,a,b,z\}$, $|I_v|=3$, and where $\{x,y\}$ and $\{a,b\}$ satisfy (\ref{potential_multihyperedges}) and are thus potential multi-hyperedges. To output a reconfiguration $I'_v,I'_w$, that does not create multi-edges, the following constraints need to hold:
\begin{align*}
    m_{I'_v}(x) &\neq m_{I'_v}(y)\\
    m_{I'_v}(a) &\neq m_{I'_v}(b).
\end{align*}
Creating the partition sequentially, i.e., by first picking a random number of copies of $x$ to go into $I'_v$, then $y, a, b$ and then $z$ (in each step making sure the constraints are met), creates bias. Depending on the number of copies of $x$ that is picked for $I'_v$, there are different numbers of partitions possible, which will be picked with a different probability:
\begin{itemize}
    \item Picking 0 copies of $x$ in $I'_v$: there are 3 possible partitions that satisfy the constraints;
    \item Picking 1 copy of $x$ in $I'_v$: there are 3 possible partitions that satisfy the constraints;
    \item Picking 2 copies of $x$ in $I'_v$: there are 2 possible partitions that satisfy the constraints.
\end{itemize}
In particular, the partitions $I'_v = \{x,x,a\}$ and $I'_v = \{x,x,b\}$ are both picked with probability 1/6, whereas the other partitions are picked with probability 1/9. 

If the partition were to be sequentially created by first picking the number of copies of $a$, then $y,x,b$ and then $z$, we obtain partition $I'_v = \{a,a,x\}$ and $I'_v = \{a,a,b\}$ with probability 1/6 and the other partitions with probability 1/9.

\subsection{\texttt{hypertrade\_nodeg} not creating multi-edges in directed hypergraphs in \texorpdfstring{$\mathcal{H}_{s}(\vb*{d}_2)$}{H-s}}
Consider the example where 
\begin{align*}
    I^{\tail}_{v} \cup I^{\tail}_w &= \{x,y,a,b,m\},\\
    I^{\head}_{v} \cup I^{\head}_w &= \{x,y,a,b,n\},
\end{align*}
and $|I^{\tail}_{v}| = 2, |I^{\head}_{v}|=3$, and $\{x,y\}$ and $\{a,b\}$ satisfy (\ref{potential_multihyperedges}). To output a reconfiguration $I'_v,I'_w$, that does not create multi-edges, the following constraints need to hold:
\begin{align*}
    (m_{I^{\prime \tail}_v}(x) &\neq m_{I^{\prime \tail}_v}(y)) \lor (m_{I^{\prime \head}_v}(x) \neq m_{I^{\prime \head}_v}(y))\\
    (m_{I^{\prime \tail}_v}(a) &\neq m_{I^{\prime \tail}_v}(b)) \lor (m_{I^{\prime \head}_v}(a) \neq m_{I^{\prime \head}_v}(b)).
\end{align*} 
Like in the previous subsection, creating the partition sequentially creates bias. To illustrate, first we pick the number of copies of $x$ that goes into $I^{\prime \tail}_v$ and $I^{\prime \head}_v$, and then $y, a, b, m$ and $n$ (such that the constraints are met at each step). There are 12 possible ways to pick copies of $x$. Picking $m_{I^{\prime \tail}_v}(x) = 1$ and $m_{I^{\prime \head}_v}(x)=0$ leaves 2 possible partitions that satisfy the constraints. Picking $m_{I^{\prime \tail}_v}(x) = 0$ and $m_{I^{\prime \head}_v}(x)=1$ leaves 8 possible partitions that satisfy the constraints. In particular, we obtain the partition $I^{\prime \tail}_v = \{x,a\}$, $I^{\prime \head}_v = \{a,b,n\}$ with probability 1/24 and the partition $I^{\prime \tail}_v = \{a,b\}$, $I^{\prime \head}_v = \{x,a,n\}$ with probability 1/96. 

\subsection{\texttt{hypertrade\_nodeg} not creating self-loops in directed hypergraphs in \texorpdfstring{$\mathcal{H}_{m}(\vb*{d}_2)$}{H-m}}
Consider the example where 
\begin{align*}
    I^{\tail}_{v} \cup I^{\tail}_w &= \{x,y,z\},\\
    I^{\head}_{v} \cup I^{\head}_w &= \{x,y,a\},
\end{align*}
and $|I^{\tail}_{v}| = 2, |I^{\head}_{v}|=1$, and $x$ and $y$ satisfy 
\begin{align*}
    &x \in I^{\tail}_v \cup I^{\tail}_w\\
    &x \in I^{\head}_v \cup I^{\head}_w\\
    &\forall u \in V \backslash \{v,w\}: x \in I^{\tail}_u \iff x \in I^{\head}_u,
\end{align*}
so that $x$ and $y$ can become self-loops after the hypertrade with nodes $u$ and $v$. To output a reconfiguration $I'_v,I'_w$, that does not create self-loops, $x$ and $y$ should not appear in the head and tail incidence set of the same node, so there must hold
\begin{align*}
    &x,y \in I^{\prime \tail}_v \cup I^{\prime \head}_v\\
    &x,y \in I^{\prime \tail}_w \cup I^{\prime \head}_w.\\
\end{align*} 
Like in the previous subsection, creating the partition sequentially creates bias. To illustrate, first we pick the number of copies of $x$ that goes into $I^{\prime \tail}_v$ and $I^{\prime \head}_v$, and then $y, z$ and $a$ (such that the constraints are met at each step). There are 2 possible ways to pick copies of $x$. Picking $x \in I^{\prime \tail}_v$ and $x \notin  I^{\prime \head}_v$ leaves 2 possible partitions that satisfy the constraints. Picking $x \notin I^{\prime \tail}_v$ and $x \in  I^{\prime \head}_v$ leaves 1 possible partition that satisfies the constraints. In particular, we obtain the partition $I^{\prime \tail}_v = \{x,y\}$, $I^{\prime \head}_v = \{a\}$ with probability 1/4 and the partition $I^{\prime \tail}_v = \{y,z\}$, $I^{\prime \head}_v = \{x\}$ with probability 1/2. 

\section{Perturbation degree}
\label{app:perturbation}
The perturbation degree of a directed hypergraph $H$ w.r.t. $G$ is 
\\
\begin{definition}[Perturbation degree for directed hypergraphs]
    Let $G,H \in \mathcal{H}_{x} (\vb*{d})$ be directed and let $d_v = (d_v^{\textnormal{out}}, d_v^{\textnormal{in}})$ for all $v \in V$. The perturbation degree of $H$ w.r.t. $G$ is
    \begin{align*}
        p_G(H) = \frac{1}{2\sum_{v \in V} (d_v^{\textnormal{out}} + d_v^{\textnormal{in}})} \sum_{v \in V} \sum_{e \in E(H)} \Big( |m_{e^{\tail}}(v) - m_{\Tilde{e}^{\tail}}(v)| + |m_{e^{\head}}(v) - m_{\Tilde{e}^{\head}}(v)| \Big).
    \end{align*}
\end{definition}

To compute the perturbation degree and approximate the mixing times of the Hypercurveball and hyperedge-shuffle method, we label the hyperedges in the initial hypergraph. Then, after a step in one of the algorithms, the modified hyperedges receive the same labels as the original hyperedges had. These labels are used as the mapping of the hyperedges in $H$ to the hyperedges in $G$.

Since we consider the hypergraphs to have unlabeled hyperedges, we are aware that the mapping that we used is not optimal. The best mapping to use for the hyperedges in $H$ to those in $G$ is the mapping that minimizes the perturbation degree. Let $\mathcal{F}$ be the set of all one-to-one mappings of hyperedges from $H$ to $G$. The perturbation degree can then be defined as
\begin{align*}
    p_G(H) = \frac{1}{2\sum_{v \in V} d_v} \min_{f \in \mathcal{F}} \sum_{v \in V} \sum_{e \in E(H)} |m_e(v) - m_{f(e)}(v)|
\end{align*}
for undirected hypergraphs, and similarly for directed hypergraphs. However, computing such optimal mapping is computationally expensive. Therefore, we choose to use the labeling system as described before as an approximation of the perturbation degree. Asymptotically, this should not matter, as both definitions of the perturbation degree converge to the same value.

\section{Minimum degree of two nodes/hyperedges}
\label{app:minimum_of_two}
We show how to compute the minimum degree of two randomly chosen nodes. A similar computation holds for the minimum degree of two hyperedges.

In a dataset, let the number of nodes with degree $d$ be given by $N(d)$. Let $\max_i d(v_i) = d_{\max}$ be the biggest node degree and let $n=|V|$ be the number of nodes. Let $v_1,v_2$ be two nodes sampled uniformly at random from the dataset (without replacement) and let $Y = \min(d_{v_1},d_{v_2}).$ Then,
\begin{align}
\label{eq:sampling_without_replacement}
    F_Y(y) &:= P(Y \leq y) = 1-P(Y > y) \nonumber \\
    &= 1-P(d_{v_1}>y)P(d_{v_2}>y|d_{v_1}>y) \nonumber\\
    &= 1-\Big(\frac{1}{n}\sum_{d_1=y+1}^{d_{\max}} N(d_1)\Big)\Big(\frac{1}{n-1}\big(\sum_{d_2=y+1}^{d_{\max}}N(d_2) - 1\big)\Big) \nonumber\\
    &= 1- \frac{1}{n(n-1)} \sum_{d_1=y+1}^{d_{\max}} N(d_1) \Big(\sum_{d_2=y+1}^{d_{\max}} N(d_2) -1 \Big).
\end{align}
Moreover,
\begin{align*}
    \mathds{E}_V[Y] = \sum_{y=0}^{d_{\max}} (1-F_Y(y)).
\end{align*}

For large datasets, $F_Y(y)$ can be approximated by considering $v_1$ and $v_2$ to be sampled \emph{with} replacement. Then,
\begin{align}
\label{eq:sampling_with_replacement}
    F_Y(y) &\approx 1-P(d_{v_1}>y)P(d_{v_2}>y) \nonumber\\
    &= 1-\Big(\frac{1}{n}\sum_{d_1=y+1}^{d_{\max}} N(d_1)\Big)^2 \nonumber\\
    &= 1-(1-F_{d_{v_1}}(y))^2,
\end{align}
where $F_{d_{v_1}}(y):=P(d_{v_1} \leq y)$.

Table \ref{tab:with_or_without_replacement} shows the values of $\mathds{E}_V[Y]$ that are obtained using either Equation (\ref{eq:sampling_without_replacement}) or Equation (\ref{eq:sampling_with_replacement}). The values show that the error of the approximation is very small. Nonetheless, the exact computation is used in this work.

\begin{table}[h]
    \centering
    \caption{Comparing the computations using Equation (\ref{eq:sampling_without_replacement}) with the approximation using Equation (\ref{eq:sampling_with_replacement}).}
    \begin{tabular}{c|c|c}
        Dataset & $\mathds{E}_V[Y]$ using Eq. (\ref{eq:sampling_without_replacement}) & $\mathds{E}_V[Y]$ using approximation Eq. (\ref{eq:sampling_with_replacement})\\ \hline
        \texttt{unicodelang} \cite{kunegis2013} & 2.34 & 2.35 \\
       \texttt{board\_directors} \cite{seierstad2011} & 1.01 & 1.01 \\
    \texttt{email\_eu} \cite{leskovec2007}  & (12.00, 33.79) & (12.03, 33.84) \\
    \texttt{thiol\_oscillator} \cite{semenov2016} & (2.75, 2.21) & (2.87, 2.35) \\
    \texttt{iAF1260b} \cite{yadati2020} & (0.97, 0.76) & (0.97, 0.77) \\ 
    \texttt{artificial\_data\_1} & 16.00 & 16.01 \\ \hline
    \end{tabular}
    \label{tab:with_or_without_replacement}
\end{table}

\newpage

\section{Artificial datasets}
\label{app:artificial_dataset}
\subsection{\texttt{artificial\_data\_1}}
The \texttt{artificial\_data\_1} dataset has node set $V = \{v_0,v_1,\hdots,v_{50}\}$ and hyperedge set $E=\{e_0,e_1,\hdots,e_{50}\}$. Nodes $v_0, v_2, \hdots, v_{49}$ have degree 16 and node $v_{50}$ has degree 30. Hyperedges $e_0, e_1, \hdots, e_{24}$ have degree 2 and hyperedges $e_{25}, e_{26}, \hdots, e_{50}$ have degree 30.

The hyperedges are:
\begin{table}[h!]
    \centering
    \begin{tabular}{ll}
        $\{2v_i\}$ & for $i \in \{0,1,2\}$ for $j \in \{0,1,\hdots, 7\}$ \\
        $\{2v_3\}$ &\\
        $\{14v_{3+i}, 16v_{7+8i}\}$ & for $i \in \{0,1,2,3\}$\\
        $\{2v_{4+i}, 16v_{14+8i}, 12v_{47+i}\}$ & for $i \in \{0,1,2\}$\\
        $\{16v_{8(j+1)+i}, (2i)v_{29+j+3i}, (14-2i)v_{32+j+3i}\}$ & for $i \in \{0,1,\hdots, 5\}$ for $j \in \{0,1,2\}$\\
        $\{30v_{50}\}$. &
    \end{tabular}
\end{table}

\subsection{\texttt{artificial\_data\_2}}
The \texttt{artificial\_data\_2} dataset has node set $V = \{v_0,v_1,\hdots,v_{1000}\}$ and hyperedge set $E=\{e_0,e_1,\hdots,e_{1000}\}$. Nodes $v_0, v_2, \hdots, v_{499}$ have degree 1 and nodes $v_{500}, v_{501}, \hdots, v_{1000}$ have degree 99. Hyperedges $e_0, e_1, \hdots, e_{999}$ have degree 50 and hyperedge $e_{1000}$ has degree 99.

The hyperedges are:
\begin{table}[h!]
    \centering
    \begin{tabular}{ll}
        $\{v_{50i}, v_{50i+1}, \hdots, v_{50i+49}\}$ & for $i \in \{0,1,\hdots,9\}$ \\
        $\{50v_{500+i}\}$ & for $i \in \{0,1,\hdots,499\}$\\
        $\{49v_{500+i+49j}, v_{999-j}\}$ & for $i \in \{0,1,\hdots,48\}$ for $j \in \{0,1,\hdots,9\}$\\
        $\{99v_{1000}\}$. & 
    \end{tabular}
\end{table}

\subsection{\texttt{artificial\_data\_3}}
The \texttt{artificial\_data\_3} dataset has node set $V = \{v_0,v_1,\hdots,v_{500}\}$ and hyperedge set $E=\{e_0,e_1,\hdots,e_{500}\}$. Nodes $v_0, v_2, \hdots, v_{249}$ have degree 1 and nodes $v_{250}, v_{251}, \hdots, v_{500}$ have degree 49. Hyperedges $e_0, e_1, \hdots, e_{499}$ have degree 25 and hyperedge $e_{500}$ has degree 49.

The hyperedges are:
\begin{table}[h!]
    \centering
    \begin{tabular}{ll}
        $\{v_{25i}, v_{25i+1},\hdots,v_{25i+24}\}$ & for $i \in \{0,1,\hdots,9 \}$\\
        $\{25v_{250+i}\}$ & for $i \in \{0,1,\hdots,249\}$\\
        $\{24v_{250+i+24j}, v_{499-j}\}$ & for $i \in \{0,1,\hdots,23\}$ for $j \in \{0,1,\hdots,9\}$\\
        $\{49v_{500}\}$. & 
    \end{tabular}
\end{table}

\section{Data sets}
\label{app:more_datasets}
The statistics of the 17 extra data sets used to validate the mixing time hypothesis are shown in Table \ref{tab:more_datasets}. For each data set, the numerical experiments performed to generate Figure \ref{fig:MP_all_data} are averaged over 100 runs, except for the data sets \texttt{NYC\_restaurant\_checkin} and \texttt{NYC\_restaurant\_tips} that are averaged over 10 runs.

\begin{table}[tbp]
\centering
\begin{threeparttable}
    \centering
    \caption{Extra data sets used for experiments on mixing time.}
    \begin{tabular}{|c c|c|c|c|c|}
    \hline
       Data set & & \# & Average degree & Median degree & $\mathds{E}[\min\{d_1,d_2\}]$\tnote{1}\\ \hline \hline
       \texttt{ceo\_club} & nodes & 25 & 3.80 & 3.00 & 3.09 \\
       \cite{data_ceo_club} & hyperedges & 15 & 6.33 & 4.00 & 4.01 \\ \hline
       \texttt{elite} & nodes & 20 & 4.95 & 5.00 & 3.75 \\
        \cite{data_elite} & hyperedges & 24 & 4.13 & 3.00 & 2.62 \\ \hline
    \texttt{kidnappings} & nodes & 246 & 1.63 & 1.00 & 1.10 \\
    \cite{data_kidnappings} & hyperedges & 105 & 3.83 & 2.00 & 1.78\\ \hline
    \texttt{plant\_pol\_kato} & nodes & 679 & 1.78 & 1.00 & 1.11 \\
    \cite{data_plant_pol_kato}& hyperedges & 91 & 13.25 & 4.00 & 4.34 \\ \hline
    \texttt{crime} & nodes & 829 & 1.78 & 1.00 & 1.10 \\ 
    \cite{data_crime} & hyperedges & 551 & 2.68 & 2.00 & 1.76 \\ \hline
    \texttt{plant\_pol\_robertson} & nodes & 456 & 33.48 & 18.00 & 13.37 \\
    \cite{data_plant_pol_robertson} & hyperedges & 1428 & 10.69 & 3.00 & 3.49\\ \hline
    \texttt{NYC\_restaurant\_checkin} & nodes & 2060 & 13.18 & 8.00 & 6.04 \\
    \cite{data_NYC_restaurant_checkin_tips}& hyperedges & 2876 & 9.44 & 6.00 & 4.55 \\ \hline
    \texttt{NYC\_restaurant\_tips} & nodes & 3112 & 3.33 & 1.00 & 1.45 \\ 
    \cite{data_NYC_restaurant_checkin_tips} & hyperedges & 3289 & 3.15 & 2.00 & 1.80 \\ \hline
    \texttt{contact\_primary\_school} & nodes & 242 & 126.98 & 123.50 & 95.37 \\
    \cite{data_contact_primary_school} & hyperedges & 12704 & 2.42 & 2.00 & 2.15\\ \hline
    \texttt{senate\_bills} & nodes & 294 & 789.62 & 592.50 & 448.54 \\
    \cite{data_senate_bills}& hyperedges & 29157 & 7.96 & 4.00 & 3.71 \\ \hline
    \texttt{house\_committees} & nodes & 1290 & 9.20 & 7.00 & 5.36 \\ 
    \cite{data_house_committees} & hyperedges & 341 & 34.79 & 40.00 & 22.57 \\ \hline
    \texttt{cat\_edge\_music\_} & nodes & 1106 & 9.49 & 5.00 & 4.64 \\
    \texttt{ blues\_reviews} \cite{data_cat_edge_music_blues_reviews} & hyperedges & 694 & 15.13 & 10.00 & 7.76\\ \hline
    \texttt{cat\_edge\_madison\_} & nodes & 565 & 8.14 & 7.00 & 5.12 \\
    \texttt{restaurant\_reviews} \cite{data_madison_vegas_algebra_geometry}& hyperedges & 601 & 7.66 & 5.00 & 4.11 \\ \hline
    \texttt{cat\_edge\_vegas\_} & nodes & 1234 & 9.62 & 8.00 & 6.96 \\ 
    \texttt{bars\_reviews} \cite{data_madison_vegas_algebra_geometry} & hyperedges & 1194 & 9.94 & 5.00 & 4.31 \\ \hline
    \texttt{cat\_edge\_} & nodes & 423 & 19.53 & 10.00 & 7.04 \\
    \texttt{algebra\_questions} \cite{data_madison_vegas_algebra_geometry} & hyperedges & 1268 & 6.52 & 4.00 & 3.76\\ \hline
    \texttt{cat\_edge\_} & nodes & 580 & 21.53 & 10.00 & 7.49 \\
    \texttt{geometry\_questions} \cite{data_madison_vegas_algebra_geometry}& hyperedges & 1193 & 10.47 & 4.00 & 4.25 \\ \hline
    \texttt{uchoice\_walmart\_items} & nodes & 183 & 140.60 & 46.00 & 50.91 \\ 
    \cite{data_uchoice_walmart_items} & hyperedges & 16698 & 1.54 & 1.00 & 1.24 \\ \hline
    \end{tabular}
     \begin{tablenotes}
    \item[1] The expected minimum degree of two randomly chosen nodes/hyperedges.
  \end{tablenotes}
\label{tab:more_datasets}
\end{threeparttable}
\end{table}

\section{Mixing times}
\label{app:fit}
\begin{figure}[tbp]
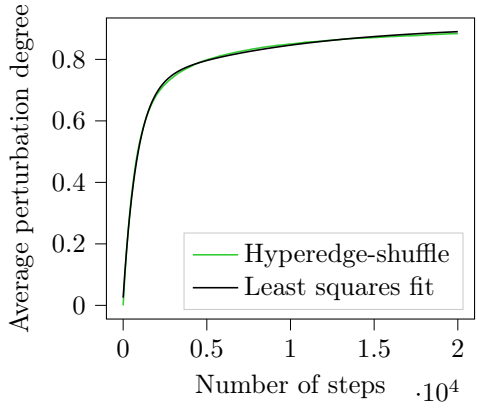
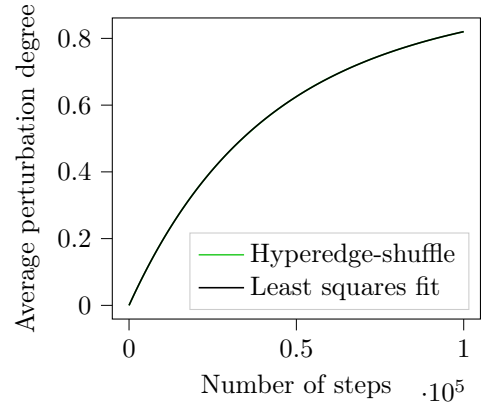
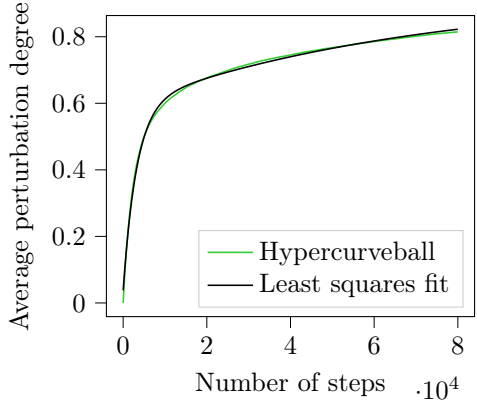
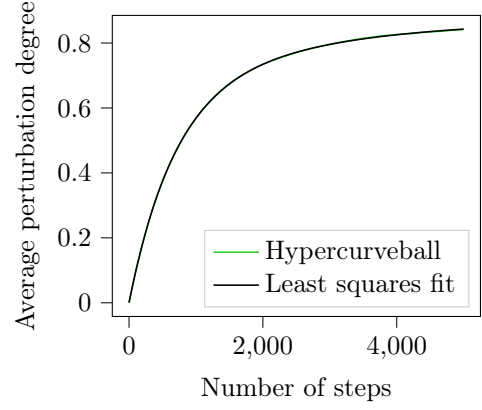

\centering
\begin{subfigure}{0.43\textwidth}
    \centering
    \subfile{Tikz/fit_unicodelang}
    \caption{\texttt{unicodelang}. $\textnormal{RMSE} = 5.01\cdot 10^{-3}$}
\end{subfigure} \hspace{0.1 \textwidth}
\begin{subfigure}{0.43\textwidth}
    \centering
    \subfile{Tikz/fit_email_eu.tex}
    \caption{\texttt{email\_eu}. $\textnormal{RMSE} = 1.93\cdot 10^{-4}$}
\end{subfigure}\\
\begin{subfigure}{0.43\textwidth}
    \centering
    \subfile{Tikz/fit_iAF1260b.tex}
    \caption{\texttt{iAF1260b}. $\textnormal{RMSE} = 7.48\cdot 10^{-3}$}
\end{subfigure} \hspace{0.1 \textwidth}
\begin{subfigure}{0.43\textwidth}
    \centering
    \subfile{Tikz/fit_plant_pol_kato.tex}
    \caption{\texttt{plant\_pol\_kato}. $\textnormal{RMSE} = 9.10\cdot 10^{-4}$}
\end{subfigure}
\caption{Results of fit $f(x)=L - a\exp(-bx)-c\exp(-dx)$ to experiments, using least squares.}
\label{fig:fits}
\end{figure}

We fit the function $f(x)=L-a\exp(-bx) -c\exp(-dx)$ to the experiments that did not show stabilizing of the perturbation degree yet, using least squares. The parameter $L$ is the limiting value of the perturbation degree, which we know for the experiments in which only one of the methods did not show a stabilized perturbation degree. For 4 extrapolated experiments, the fitting parameters are given in Table \ref{tab:fitting_parameters}. A visualization of these fits and their Root Mean Square Error (RMSE) is shown in Figure \ref{fig:fits}.
\begin{table}[h]
    \centering
    \caption{Fitting $f(x)=L-a\exp(-bx) -c\exp(-dx)$ using least squares.}
    \begin{tabular}{c|c|c|c|c|c|c}
        Data set & Method & $a$ & $b$ & $c$ & $d$ & $L$\\ \hline
        \texttt{unicodelang} & Hyperedge-shuffle & 1.98e-1 & 1.10e-4 & 6.89e-1 & 1.20e-3 & 9.12e-1 \\ 
        \texttt{email\_eu} & Hyperedge-shuffle & 3.09e-1 & 8.77e-6 & 6.79e-1 & 2.85e-5 & 9.88e-1\\ 
        \texttt{iAF1260b} & Hypercurveball & 3.35e-1 & 1.46e-5 & 5.54e-1 & 3.10e-4 & 9.27e-1\\
        \texttt{plant\_pol\_kato} & Hypercurveball & 2.30e-1 & 3.75e-4 & 6.49e-1 & 1.45e-3 & 8.79e-1
    \end{tabular}
    \label{tab:fitting_parameters}
\end{table}

The estimated mixing times of all experiments are shown in Table \ref{tab:mixing_times}.

\begin{table}[tbp]
\centering
\begin{threeparttable}
    \centering
    \caption{Mixing times of all datasets.}
    \begin{tabular}{c|c|c}
        Data set & Mixing time Hypercurveball & Mixing time hyperedge-shuffle \\ \hline
        \texttt{unicodelang} & 3136 & 21665\tnote{1}\\
        \texttt{board\_directors} & 4519 & 970 \\
        \texttt{email\_eu} & 12027 &314051\tnote{1} \\
        \texttt{thiol\_oscillator} & 103 & 229\\
        \texttt{iAF1260b} & 198835\tnote{1}& 9696\\
        \texttt{artificial\_data\_1} & 201 & 328\\
        \texttt{artificial\_data\_2} & 7575 & 3915\\
        \texttt{artificial\_data\_3} & 3736 & 1952\\
         \texttt{ceo\_club} & 87 & 58 \\
        \texttt{elite} & 81 &110 \\
     \texttt{kidnappings} & 1988\tnote{1}& 1673 \\
     \texttt{plant\_pol\_kato} & 6860\tnote{1}& 1459 \\
     \texttt{crime} & 8683 & 3446 \\ 
     \texttt{plant\_pol\_robertson} & 4630 & 20404\tnote{1}\\
     \texttt{NYC\_restaurant\_checkin} & 23442 & 32220 \\
     \texttt{NYC\_restaurant\_tips} & 77877\tnote{1}& 24996 \\ 
     \texttt{contact\_primary\_school} & 1118 & 63781\tnote{1}  \\
     \texttt{senate\_bills} & 1806 & 178073\tnote{1}  \\
     \texttt{house\_committees} & 7881 & 1782 \\ 
     \texttt{cat\_edge\_music\_blues\_reviews} & 9329 & 5247 \\
     \texttt{cat\_edge\_madison\_restaurant\_reviews} & 4107 & 4418 \\
     \texttt{cat\_edge\_vegas\_bars\_reviews} & 7204 & 13879 \\ 
     \texttt{cat\_edge\_algebra\_questions} & 6264 & 9065 \\
     \texttt{cat\_edge\_geometry\_questions} & 7939 & 13141 \\
     \texttt{uchoice\_walmart\_items} & 2865 & 78734 \\ 
    \end{tabular}
     \begin{tablenotes}
    \item[1] Value extrapolated.
  \end{tablenotes}
\label{tab:mixing_times}
\end{threeparttable}
\end{table}
\newpage
\printbibliography

\end{document}